\documentclass[a4,10pt,twocolumn]{scrartcl}
\usepackage[top=2.4cm, bottom=2.8cm, left=1.6cm, right=1.6cm]{geometry}
\pdfoutput=1

\input{myStandardPreamble_arXiv.tex}
\usepackage[english]{babel}
\usepackage[utf8x]{inputenc}
\usepackage{enumitem}					

\usepackage{libertine}
\usepackage[T1]{fontenc}
\usepackage[noBBpl]{mathpazo}

\usepackage[hang]{footmisc}
\setfnsymbol{wiley}

\makeatletter
\renewcommand{\@biblabel}[1]{[#1]\hfill}
\makeatother

\makeatletter
\let\NAT@parse\undefined
\makeatother

\usepackage[format=plain,			
            labelsep=period,		
            font=small,			
            labelfont=bf,			
            skip=5pt			
            ]{caption}
            
\usepackage{etoolbox}
\apptocmd{\thebibliography}{\setlength{\itemsep}{2pt}}{}{}

\usepackage[final,kerning=false,protrusion=false]{microtype}

\usepackage{authblk}

\newtheorem{theorem}{Theorem}
\newtheorem{lemma}{Lemma}
\newtheorem{definition}{Definition}
\newtheorem{remark}{Remark}
\newtheorem{example}{Example}
\newtheorem{assumption}{Assumption}
\newtheorem{corollary}{Corollary}
\newtheorem{proposition}{Proposition}

\let\oldTheorem\theorem
\renewcommand{\theorem}{\oldTheorem\normalfont}
\let\oldLemma\lemma
\renewcommand{\lemma}{\oldLemma\normalfont}
\let\oldCorollary\corollary
\renewcommand{\corollary}{\oldCorollary\normalfont}
\let\oldDefinition\definition
\renewcommand{\definition}{\oldDefinition\normalfont}
\let\oldRemark\remark
\renewcommand{\remark}{\oldRemark\normalfont}
\let\oldAssumption\assumption
\renewcommand{\assumption}{\oldAssumption\normalfont}
\let\oldExample\example
\renewcommand{\example}{\oldExample\normalfont}

\addtokomafont{paragraph}{\itshape}

\colorlet{fictitiousColor}{green!45!black}    

\newboolean{showCommentsOld}
\setboolean{showCommentsOld}{false}

\newcommand{\Real}{\text{Re}}								
\newcommand{\unitVec}[2][{}]{{}_{#1} e_{#2}}                		
\newcommand{\img}{\mathrm{i} }								
\newcommand{\rank}{\textup{rank}}							
\newcommand{\map}[3]{#1:#2 \rightarrow #3}						
\newcommand{\ones}{\mathbf{1}}								
\newcommand{\nat}{\mathbb{N}}								
\newcommand{\natPos}{\mathbb{N}_{>0}}							
\newcommand{\real}{\mathbb{R}}								
\newcommand{\realNonNeg}{\mathbb{R}_{\geq 0}}					
\newcommand{\realPos}{\mathbb{R}_{>0}}						
\newcommand{\C}{\mathcal{C}}								
\newcommand{\pder}[2]{\frac{\partial #1}{\partial #2}}			
\newcommand{\sign}{\textup{sgn}}							
\newcommand{\oprocendsymbol}{\hbox{$\bullet$}} 					
\newcommand{\oprocend}{\relax\ifmmode\else\unskip\hfill\fi\oprocendsymbol}

\newcommand{\ineqConstrL}{C}								
\newcommand{\ineqConstrLAgent}{\MakeLowercase{\ineqConstrL}}		%
\newcommand{\ineqConstrR}{d}		                        		
\newcommand{\eqConstrL}{A}			                        		
\newcommand{\eqConstrLAgent}{\MakeLowercase{\eqConstrL}}    		%
\newcommand{\eqConstrR}{b}			                        		
\newcommand{\setEq}{\mathcal{I}_{\textup{eq}}}					
\newcommand{\setIneq}{\mathcal{I}_{\textup{ineq}}}				
\newcommand{\lagrangian}{L}		                            	
\newcommand{\dualIneq}{\lambda}                            		
\newcommand{\dualEq}{\nu}                                   		

\newcommand{\laplacian}{G}			                        		
\newcommand{\elLaplacian}{\MakeLowercase{\laplacian}}			
\newcommand{\G}{\mathcal{G}}								
\newcommand{\Adj}{\mathbf{A}}								
\newcommand{\AdjEl}{\mathbf{\MakeLowercase{\Adj}}}				
\newcommand{\degreeMat}{\mathrm{D}}							
\newcommand{\degreeMatEl}{\mathrm{\MakeLowercase{\degreeMat}}}		
\newcommand{\subpath}{\mathrm{subpath}}                     		
\newcommand{\pathLeft}{\langle}								
\newcommand{\pathRight}{\rangle}							
\newcommand{\pathSep}{|}									
\newcommand{\length}{\ell}									
\newcommand{\tail}{\textup{tail}}							
\newcommand{\head}{\textup{head}}							

\newcommand{\degree}{\delta}                                		
\newcommand{\lFac}{\textup{left}}                           		
\newcommand{\rFac}{\textup{right}}                          		
\newcommand{\resort}{\textup{proj}}							
\newcommand{\fBr}{\mathcal{FB}r}							
\newcommand{\lieBr}{\mathcal{LB}r}							
\newcommand{\Ev}{\textup{Ev}}								

\newcommand{\state}{z}									
\newcommand{\seqParam}{\sigma}								
\newcommand{\recBracket}[2]{R_{#1}}							
\newcommand{\recBracketPHall}[2]{\tilde{R}_{#1}}				

\newcounter{MyAssumptionCounter}
\renewcommand{\emph}{\textit}

\ifthenelse{\boolean{showCommentsOld}}{\colorlet{commentColor}{orange}}{\colorlet{commentColor}{black}}
\ifthenelse{\boolean{showCommentsOld}}{}{\renewcommand{\sout}{}}

\newcommand{\simonResub}[1]{{\color{commentColor} #1}}

\newcommand{\margin}[1]{\marginpar{\color{red}\tiny\ttfamily#1}}

\renewcommand{\margin}[1]{\marginpar{\color{red}\tiny\ttfamily}}

\newboolean{longVersion}
\setboolean{longVersion}{true}

\title{\LARGE \textbf
A Lie bracket approximation approach to distributed optimization over directed graphs
\footnote{This article is a sligthly extended version of \cite{mic2017distributed} with the follwing additions: an extra illustration in~\Cref{figIllustrationSubpaths}, an additional result~\Cref{lemmaEquivalentBrackets}\simonResub{, some supplementary details on formal brackets in~\Cref{secAppendixFormalBrackets}, a section on filtered sadle-point dynamics (\Cref{secFilteredSPD}), a slightly extended example which illustrates the effect of such filtered dynamics and a proof of~\Cref{lemmaSaddlePointDynamics}.}}
}

\setkomafont{section}{\Large}
\setkomafont{subsection}{\large}
\setkomafont{subsubsection}{\itshape}

\begin{document}

\date{}
\author[1]{Simon Michalowsky}
\author[2]{Bahman Gharesifard}
\author[1]{Christian Ebenbauer}
\affil[1]{Institute for Systems Theory and Automatic Control, University of Stuttgart, Germany \protect\\ \texttt{\small $\lbrace$michalowsky,ce$\rbrace$@ist.uni-stuttgart.de}}
\affil[2]{Department of Mathematics and Statistics, Queen's University, Canada \protect\\ \texttt{\small bahman@mast.queensu.ca} \protect\\[2em]}

\maketitle

\begin{abstract}
\textbf{Abstract.} We consider a group of computation units trying to cooperatively
	solve a distributed optimization problem with shared linear
	equality and inequality constraints. Assuming that the 
	computation units are communicating over a network whose
	topology is described by a time-invariant directed graph,
	by combining saddle-point dynamics with Lie bracket approximation techniques
	we derive a methodology that allows to design distributed 
	continuous-time optimization algorithms that solve this problem
	under minimal assumptions on the graph topology as well as on
	the structure of the constraints. We discuss several extensions
	as well as special cases in which the proposed procedure
	becomes particularly simple.
\end{abstract}

\section{Introduction}
Driven by new applications and advancing communication technologies,
the idea of solving optimization problems in a distributed fashion
using a group of agents interchanging information over a communication
network has gained a lot of interest during the last decades. 
Application examples include, among others, optimal power dispatch
problems in {smart grids} \cite{geidl2005optimal}, distributed machine learning \cite{boyd2011distributed} or formation
control problems \cite{bullo2009robotic}. 
{Besides several results on distributed computation~\cite{costello2014degree}, controllability and stabilization~\cite{belabbas2013sparse, chen2015controllability, bahman2017stabilization},}
there also exists a vast body of literature on distributed optimization
algorithms, both in discrete- \cite{nedic2015distributed, boyd2011distributed} {and} continuous-time 
\cite{feijer2010stability, wang2011control, duerr2012saddlePoint, niederlaender2015distributed, gharesifard2014distributed, touri2016saddle}, where in the
present work we will focus on the latter one. While in most of
the works a consensus-based approach is used where all agents aim to agree 
on a common solution of the overall optimization problem, in the last
years other solutions have been proposed as well \cite{niederlaender2015distributed}.
However, it is usually assumed that the underlying communication 
network is of undirected nature {or is weight-balanced} and it has turned out that
establishing distributed optimization algorithms in the presence of
directed communication structures is much more difficult.
While there exist some approaches aiming to address this problem \cite{gharesifard2014distributed, touri2016saddle},
these are limited to unconstrained optimization problems using
a consensus-based approach.

{The contribution of this work is to provide} {a unified framework that allows} {the design of} continuous-time
distributed optimization {algorithms} {for} a very general class of constrained
optimization problems under mild assumptions on the possibly
directed underlying communication network. 
The main idea of our
approach is to employ classical saddle-point dynamics with
proven convergence guarantees in a centralized setting and
derive distributed approximations thereof. To this end,
we follow a {two-step} procedure where we first propose suitable 
Lie bracket representations of saddle-point dynamics and 
then use ideas from geometric control theory to design distributed 
approximations thereof. This idea has already been employed
in previous works using a consensus-based approach \cite{ebenbauer2017directed}
and for more general optimization problems with linear equality constraints
in a gradient-free setting \cite{mic2017es}. However, the focus
in both works was on the first step of rewriting the saddle-point
dynamics and the second step of designing distributed approximations
was rarely treated. {In the present paper we further contribute
to both steps: on the one hand, we extend the class of optimization problems
the approach is applicable to, and, on the other, we present an 
algorithm for {designing} suitable approximations.} 
While we limit ourselves to convex
optimization problems with linear equality and inequality
constraints, we emphasize that
the same techniques may be used for a much larger class
of optimization problems{, see~\cite{mic2018extensions}.}
{We further emphasize that the main goal of this work
is not to provide distributed algorithms ready to be implemented in practical
applications but to present a general framework that in principle
allows addressing several limitations common in distributed
optimization and control problems.}

\section{Preliminaries}
\paragraph*{Notation.}
{We let $ \nat $ denote the set of non-negative integers
and let $ \natPos $ be the set of positive integers.}
{Similarly,} we denote by $ \real^n $ the set of $n$-dimensional real
vectors, by $ {\realNonNeg^n} $ those with {non-negative} entries
and by $ {\realPos^n} $ those with {positive} entries.
We further write $ \C^p $, $ p \in \nat $, for 
the set of $p$-times continuously differentiable {real-valued} functions.
The gradient of a function $ f: \real^n \to \real $, $ f \in \C^1 ${\sout{, 
with respect to its argument $x\in \real^n$,}}
will be denoted by {$ \map{\nabla f}{\real^n}{\real^n} $}{\sout{;
we often omit the subscript, if it is clear from the context}}.
We denote the $ (i,j) $th entry of a matrix $ A \in \real^{n \times m} $
by $ a_{ij} $, and sometimes denote $ A $ by $ A = [a_{ij}] $. The rank 
of $A$ is denoted by $ \rank(A) $. We use $ \unitVec{i} $ to denote the {real} vector 
with the $ i $th entry equal to $ 1 $ and all other entries equal to~$0$, {where the dimension should be clear from the context,}
and also use the short-hand notation $ \ones_n=[1,\ldots,1]^T {\in \mathbb{R}^n} $.
For a vector $ \lambda \in \real^n $ we let $ \textup{diag}(\lambda) \in \real^{n \times n} $
denote the diagonal matrix whose diagonal entries are the entries of $ \lambda $.
We denote the sign function by $ \sign: \real \to \lbrace -1,0,1 \rbrace $, 
where $ \sign(-a) = -1 $, $ \sign(a) = 1 $ for any \mbox{$ a > 0 $}
and $ \sign(0) = 0 $.
For a vector $ x=[x_1,\ldots,x_n]^T \in \real^n $ and a finite set $ S \subset \{1,\ldots, n\} $, we denote by $ x_S $ the set of all $ x_i $ with $ i\in S$. We also denote the complement of a set $ S\subset \real^n $ by $ S^c $.

Given two continuously differentiable vector fields 
$ \map{\phi_1}{\real^n}{\real^n} $ and $ \map{\phi_2}{\real^n}{\real^n} $,
the Lie bracket of $ \phi_1 $ and $ \phi_2 $ evaluated at $x$ is defined to be
\begin{align}
	[ \phi_1, \phi_2 ](x) := \pder{\phi_2}{x}(x) \phi_1(x) - \pder{\phi_1}{x}(x)  \phi_2(x). \label{eqDefLieBracket}
\end{align}
{Observe that the Lie bracket is a bilinear skew-symmetric operator
that fulfills the Jacobi-identity, see also~\cite{bourbaki1998lie}.}
For a set of vector fields $ \Phi = \lbrace \phi_1,\phi_2,\dots,\phi_M \rbrace $,
$ \phi_i: \real^n \to \real^n $, $ \phi_i \in \C^1 $,
we denote by $ \lieBr(\Phi) $ the set of Lie brackets generated
by $ \Phi $. For an (iterated) Lie bracket $ B = [B_1,B_2] $, $ B_1, B_2 \in \lieBr(\Phi) $,
we then let $ \lFac(B) = B_1 $, $ \rFac(B) = B_2 $ denote the left and
right factor of $ B $, respectively. {We note that the left and right 
factor are not uniquely defined for Lie brackets since one Lie bracket can have
multiple representations; in fact, {to obtain uniqueness,} we would need to define these operators on
the set of formal brackets of indeterminates. 
The interested reader is referred to
\ifthenelse{\boolean{longVersion}}{
\Cref{secAppendixFormalBrackets}
}
{the extended version~\cite{mic2017distributedArxiv}}
or a standard textbook such as~\cite{bourbaki1998lie}
for some more details on this subject.
}
In the following we accept this abuse of notation to avoid the formal overhead and assume that, whenever $ \lFac(B) $, $ \rFac(B) $ are used
for Lie brackets $B \in \lieBr(\Phi) $, the bracket $B$ has to interpreted 
as a formal bracket, and we assume the formal bracket representation to be given.
As an example, for the left and right factor we distinguish
between the two brackets $ \big[ \phi_1, [\phi_1,\phi_2] \big] $ and $ \big[ [\phi_2,\phi_1], \phi_1 \big] $
which are equivalent as brackets in $ \lieBr(\Phi) $ but not equivalent as
formal brackets where each bracket is a word consisting of the symbols
$ \phi_1, \phi_2 $, the brackets, as well as the comma.
We further define the degree of a Lie bracket $B \in {\mathcal{LB}}(\Phi) $ as $ \degree(B) = \tilde{\degree}_{\Phi}{(B)} $
and the degree of the $k$th vector field{, $ k = 1,2,\dots,M $,} as $ \degree_k(B) = \tilde{\degree}_{\lbrace \phi_k \rbrace}{(B)} $,
where
\begin{align*}
	\tilde{\degree}_{\mathcal{S}}(B)
	=
	\begin{cases}
		1 & \text{if } B \in \mathcal{S} \\
		{\tilde{\degree}_{\mathcal{S}}}(\lFac(B)) + {\tilde{\degree}_{\mathcal{S}}}(\rFac(B)) & \text{otherwise,}
	\end{cases}
\end{align*}
with $ \mathcal{S} \subseteq \Phi $. {Again, we note that 
formally we would require to define the degree on the set of formal
brackets for it to be mathematically precise.}

\paragraph*{Basics on graph theory.}
We recall some basic notions on graph theory, and refer the reader 
to~\cite{biggs1993algebraic} or other standard references for more information.
A directed graph (or simply digraph) is an ordered pair $ \G = ( \mathcal{V}, \mathcal{E} ) $,
where $ \mathcal{V} = \lbrace v_1, v_2, \dots, v_n \rbrace $, {$ v_i \neq v_j $ for $ i \neq j $,}
is the set of nodes and $ \mathcal{E} \subseteq \mathcal{V} \times \mathcal{V} $
is the set of edges, i.e. $ (v_i,v_j) \in \mathcal{E} $ if there 
is an edge from node $ v_i $ to $ v_j $. In our setup the
edges encode to which other agents some agent has access to, i.e. 
$ (v_i, v_j) \in \mathcal{E} $ means that node $ v_i $
receives information from node $ v_j $. We say that node $ v_j $ is an out-neighbor of
node $v_i$ if there is an edge from node $ v_i $ to node $ v_j $. 
The adjacency matrix $ \Adj = [ \AdjEl_{ij} ]  \in \real^{n\times n} $ associated to $ \G $ is defined as {
\begin{align}
	\AdjEl_{ij}
	&=
	\begin{cases}
		1 & \text{if } i \neq j \text{ and } (v_i,v_j) \in \mathcal{E},  \\
		0 \hspace*{1.3cm} & \text{otherwise}.
	\end{cases} \label{eqPreliminariesDefAdjacencyMat}
\end{align}
We also define the out-degree matrix  $ \degreeMat = [ \degreeMatEl_{ij} ] $ associated to $ \G $ as
\begin{align}
	\degreeMatEl_{ij}
	&=
	\begin{cases}
		\sum_{k=1}^{n} \AdjEl_{ik} & \text{if } i = j \\
		0 \hspace*{1.3cm}                   & \text{otherwise.}
	\end{cases}
\end{align}
Finally, we call $ \laplacian = \degreeMat - \Adj = [ \elLaplacian_{ij} ] \in \real^{n \times n} $ the Laplacian of $ \G $.}
A digraph is said to be undirected if $ (v_i,v_j) \in \mathcal{E} $
implies that $ (v_j,v_i) \in \mathcal{E} $, or, equivalently,
if $ \laplacian = \laplacian^{\top} $. Further, a digraph 
$ \G $ is called weight-balanced if $ \ones_n^T \laplacian =0 $.
A directed path in $ \G $ is a sequence of nodes connected by edges,
and we write $ p_{i_1 i_r} = \pathLeft v_{i_1} \pathSep v_{i_2} \pathSep \dots \pathSep v_{i_r} \pathRight $
for a path from node $ v_{i_1} $ to node $ v_{i_r} $. We further
denote by $ \head(p_{i_1 i_r}) = i_1$ and $ \tail(p_{i_1 i_r}) = i_r $
the {head and the tail of a path $ p_{i_1 i_r} $, respectively}.
{We also let $ \length(p_{i_1 i_r}) = r-1 $ denote the length of the path.}
A digraph $ \G $ is said to be strongly connected (or simply connected in case of
undirected graphs) if there is a directed path between any two nodes. 
For a path $p_{ij}$ from node $ v_i$ to node $v_j$ we denote by 
$ \subpath_{i \bullet}(p_{ij}) $ and $ \subpath_{\bullet j}(p_{ij}) $ 
the set of all subpaths of $ p_{ij} $ (not including $ p_{ij} $ itself) which, 
respectively, {start at $ v_i $ or end at $ v_j $}. Given a subpath 
$ q \in \subpath_{i \bullet}(p_{ij}) $, we denote by $ q^{c} $ the path in 
$ \subpath_{\bullet j}(p_{ij}) $ whose composition with $ q $ gives $ p_{ij} $.

\section{Problem setup}
Consider an optimization problem of the form
\begin{align}
	\begin{split}
		\min\limits_{x} \quad & F(x) = \sum\limits_{i=1}^n F_i(x_i) \\
		\text{s.t}      \quad & \eqConstrLAgent_i x - \eqConstrR_i = 0, \qquad i \in \setEq \subseteq \lbrace 1,2,\dots,n \rbrace, \\
		                      & \ineqConstrLAgent_i x - \ineqConstrR_i \leq 0, \qquad i \in \setIneq \subseteq \lbrace 1,2,\dots,n \rbrace, 
	\end{split}
	\label{eqOptimizationProblem}
\end{align}
where $ x = [ x_1, \dots, x_n ]^\top \in \real^n $, 
{$ \eqConstrLAgent_i, \ineqConstrLAgent_i \in \real^{1 \times n} $,}
{$ \eqConstrR_i, \ineqConstrR_i \in \real $,}
and the $ F_i: \real \to \real $, $ F_i \in \mathcal{C}^2 $,
are assumed to be strictly convex functions.
We assume further that the feasible set of \eqref{eqOptimizationProblem} is 
non-empty{; thus,} there exists a unique solution $x^* \in \real^n $
to \eqref{eqOptimizationProblem}.

The problem can be interpreted as having $n$ computation units 
or {agents} available, each one trying to optimize 
its own objective function {$ F_{i} $} while, if $ i \in \setIneq $ or $ i \in \setEq $, respecting the $i$th
global constraints among all agents. 
It is reasonable to assume that the constraints are
associated to the agents in such a way that the constraint
corresponding to agent $i$ involves its own state. This
is ensured by the following assumption on the set of constraints:
\begin{assumption}\label{assConstraints1}
{For each $ i \in \setEq $, if $ \eqConstrLAgent_i \neq 0 $,
	then $ \eqConstrLAgent_i \unitVec{i} \neq 0 $; and, for each $ i \in \setIneq $, if
	 $ \ineqConstrLAgent_i \neq 0 $,
	 then $ \ineqConstrLAgent_{i} \unitVec{i} \neq 0 $.} \oprocend
\end{assumption}
It should be noted that, {merely} for the ease of presentation, we limit
ourselves to the case that each
agent has at most one equality and one inequality constraint but the following
results apply with some modifications to the case
where each agent has several constraints, i.e., 
$ \eqConstrLAgent_i \in \real^{M_i \times n} $,
$ \ineqConstrLAgent_i \in \real^{m_i \times n} $
for some $m_i,M_i \in {\natPos} $.
{Our intention is to focus on presenting our results
in a more understandable fashion and avoid complicated notations 
introduced when considering more general problem setups.}
{Still, we emphasize that the framework is applicable
in fairly general situations, and we refer the reader to~\cite{mic2018extensions},
where we focus on the discussion of the class of distributed optimization 
problems the methodology can in principle be applied to.}

{Going along that direction of a simpler notation}, we augment the problem 
\eqref{eqOptimizationProblem} by non-restrictive constraints
such that exactly one equality and one inequality constraint is
associated to each agent, i.e., we consider the augmented problem
\begin{align}
	\begin{split}
		\min\limits_{x} \qquad & F(x) = \sum\limits_{i=1}^n F_i(x_i) \\
		\text{s.t}      \qquad & \eqConstrLAgent_i x - \eqConstrR_i = 0, \qquad i = 1,2,\dots,n, \\
		                       & \ineqConstrLAgent_i x - \ineqConstrR_i \leq 0, \qquad i = 1,2,\dots,n,
	\end{split}
	\label{eqOptimizationProblemAugmented}	
\end{align}
where $ \eqConstrLAgent_i = 0, \eqConstrR_i = 0 $ for $ i \notin \setEq $ 
and $ \ineqConstrLAgent_i = 0, \ineqConstrR_i > 0 $ for $ i \notin \setIneq $,
such that the feasible set as well as the solution of \eqref{eqOptimizationProblem}
and \eqref{eqOptimizationProblemAugmented} are the same.

{In the following, we wish to} design continuous-time algorithms 
{that ``converge''} to an arbitrarily small neighborhood
of the solution of \eqref{eqOptimizationProblemAugmented}
and that can be implemented in a \emph{distributed} 
{fashion, i.e., each agent only uses information of its} own state and objective function $F_{i}$ as
well as those of its out-neighbors, where out-neighboring agents are
defined by {a} communication graph. 

More precisely, we assume that the communication topology is given
by some directed graph $ \mathcal{G} = ( \mathcal{V}, \mathcal{E} ) $,
where $ \mathcal{V} = \lbrace v_1, v_2, \dots, v_n \rbrace $ is a finite set of nodes and
$ \mathcal{E} \subseteq \mathcal{V} \times \mathcal{V} $ 
is the set of edges between the nodes. In our setup,
the nodes play the role of the $n$ agents
and the edges define the allowed communication links 
between the agents, i.e., if there exists an edge from
agent $i$ to agent $j$, then agent $i$ has access to the
state of agent $j$. Using the graph Laplacian $ \laplacian = [ \elLaplacian_{ij} ] $
associated to $ \mathcal{\laplacian} $, we then have
the following definition of a distributed algorithm:
\begin{definition}\label{defDistributedDynamics}
We say that a continuous-time algorithm with agent dynamics
	of the form
	\begin{align}
		\dot{\state}_j = f_j(t,\state),
	\end{align}
	$ j = 1,2,\dots,{N} $, $ \state = [ \state_1, \state_2, \dots, \state_{{N}} ]^\top \in \real^{{N}} $,
	$ f_j: \real \times {\real^N} \to \real $,
	is {distributed w.r.t. the graph $ \mathcal{G} $} if 
	it can equivalently be written as
	\begin{align}
	    \dot{\state}_j = \tilde{f}_j(t,\state_{{\mathcal{N}(i)}}),
	\end{align}
	where $ \mathcal{N}(i) := \lbrace j=1,2,\dots,{N}: g_{ij} \neq 0 \rbrace $
	is the set of indices of all out-neighboring agents. \oprocend
\end{definition}
In words, $ f_i $ may only depend on $ \state_i $ and all states 
$ \state_j $ whose corresponding agent $j$ have a communication link 
to agent $j$, i.e., the algorithm obeys the
communication topology defined by the directed graph
$ \mathcal{G} $. 

Our approach relies on
the use of saddle-point dynamics, i.e., algorithms that
utilize the saddle{-}point property of the Lagrangian. The Lagrangian
$ \lagrangian: \real^n \times \real^{n} \times {\realNonNeg^{n}} \to \real $
associated to \eqref{eqOptimizationProblemAugmented} is given by 
\begin{align}
	\lagrangian(x,\dualEq,\dualIneq) 
	\ifthenelse{\boolean{longVersion}}{
	&= \sum\limits_{i=1}^{n} \big( F_{i}(x_i) +  {\dualEq_{i}}( \eqConstrLAgent_{i} x - \eqConstrR_{i}) + {\dualIneq_{i}}( \ineqConstrLAgent_{i} x - \ineqConstrR_{i} ) \big) \nonumber \\}{}
	&= F(x) + \dualEq^\top ( \eqConstrL x - \eqConstrR ) + \dualIneq^\top ( \ineqConstrL x - \ineqConstrR ), \label{eqDefLagrangian}
\end{align}
where we {have} used the stacked matrices 
\begin{align}
	\begin{array}{llllll}
	\ineqConstrL &= \begin{bmatrix} {\ineqConstrLAgent_{1}}^\top & \dots & {\ineqConstrLAgent_{n}}^\top \end{bmatrix}^\top, \quad &
	\ineqConstrR &= \begin{bmatrix} {\ineqConstrR_{1}} & \dots & {\ineqConstrR_{n}} \end{bmatrix}^\top, \\
	\eqConstrL &= \begin{bmatrix} {\eqConstrLAgent_{1}}^\top & \dots & {\eqConstrLAgent_{n}}^\top \end{bmatrix}^\top, &
	\eqConstrR &= \begin{bmatrix} {\eqConstrR_{1}} & \dots & {\eqConstrR_{n}} \end{bmatrix}^\top, \\
	\dualIneq &= \begin{bmatrix} {\dualIneq_{1}} & \dots & {\dualIneq_{n}} \end{bmatrix}^\top, &
	\dualEq &= \begin{bmatrix} {\dualEq_{1}} & \dots & {\dualEq_{n}} \end{bmatrix}^\top,
	\end{array}
\end{align}
with $ \dualEq \in \real^{n},\dualIneq \in \real^{n} $
being the associated Lagrange multipliers.
Here, a point $ (x^\star,\dualEq^\star,\dualIneq^\star) \in \real^{n} \times \real^{n} \times {\realNonNeg^{n}} $
is said to be a (global) saddle point of $L$ if 
for all $ x \in \real^n $, $ \dualEq \in \real^{n} $,
$ \dualIneq \in {\realNonNeg^{n}} $ we have
\begin{align}
    \lagrangian(x^\star,\dualEq,\dualIneq) \leq \lagrangian(x^\star,\dualEq^\star,\dualIneq^\star) \leq \lagrangian(x,\dualEq^\star,\dualIneq^\star) .
    \label{eqSaddlePointProperty}
\end{align}
It is well-known that if the Lagrangian has some saddle point
$ (x^\star,\dualEq^\star,\dualIneq^\star) $, then $ x^\star $ is a solution of 
\eqref{eqOptimizationProblemAugmented}. In the present setup, since \eqref{eqOptimizationProblemAugmented}
is a convex problem and the feasible set is non-empty, 
the existence of a saddle point is ensured (cf., e.g., \cite{hiriart2013convex}) such that
finding a saddle point of $\lagrangian$ is equivalent 
to finding a solution to \eqref{eqOptimizationProblemAugmented}.
We further require the following regularity assumption to hold:
\begin{assumption}\label{assMFCQ}
The constraints in \eqref{eqOptimizationProblem} fulfill the
	Mangasarian-Fromovitz constraint qualifications at
	the optimal solution $ x^\star $, i.e.,
	the vectors $ \eqConstrLAgent_i $, $ i \in \setEq $,
	are linearly independent and there exists $ q \in \real^n $ such 
	that $ \ineqConstrLAgent_i q < 0 $ for all $ i \in \setIneq $
	for which $ \ineqConstrLAgent_i x^\star - \ineqConstrR_i = 0 $ and
	$ \eqConstrLAgent_i q = 0 $ for all $ i \in \setEq $.
	\oprocend
\end{assumption}
This assumption ensures that the set of saddle points of the Lagrangian
associated to \eqref{eqOptimizationProblem} is non-empty and compact,
see \cite[Theorem 1]{wachsmuth2013licq}. Note that, due to the 
augmentation of the optimization problem, the set of saddle points
of the Lagrangian $ \lagrangian $ associated to \eqref{eqOptimizationProblemAugmented}
is in general not compact{, an issue that we address} by modifying the saddle-point dynamics.
To be more precise, in the following Lemma we propose a modified saddle-point dynamics,
which is an extension of the one proposed in \cite{duerr2012saddlePoint},
and show asymptotic stability of a compact subset of the set of saddle points{;
\ifthenelse{\boolean{longVersion}}{a proof is presented in \Cref{proofLemmaSaddlePointDynamics}.}{{a proof is provided in the extended version~\cite{mic2017distributedArxiv}.}}}
\begin{lemma}\label{lemmaSaddlePointDynamics}
Consider the following modified saddle-point dynamics
	\begin{subequations}
	\begin{flalign}
		\dot{x}         &= - \nabla_x L(x,\dualEq,\dualIneq)                                  &&\hspace*{-2.1em}= -\nabla F(x) \hspace*{-0.1em} -  \hspace*{-0.1em} \eqConstrL^\top  \hspace*{-0.1em} \dualEq  \hspace*{-0.1em} -  \hspace*{-0.1em} \ineqConstrL^\top  \hspace*{-0.1em} \dualIneq\hspace*{-3em} &   \label{eqSPAa} \\
		\dot{\dualEq}   &= \nabla_{\dualEq} L(x,\dualEq,\dualIneq) + w(\dualEq)               &&\hspace*{-2.1em}= \eqConstrL x - \eqConstrR + w(\dualEq) & \label{eqSPAb} \\
		\dot{\dualIneq} &= \textup{diag}(\dualIneq) \nabla_{\dualIneq} L(x,\dualEq,\dualIneq) &&\hspace*{-2.1em}= \textup{diag}(\dualIneq) (\ineqConstrL x-\ineqConstrR), &  \label{eqSPAc}
	\end{flalign}\label{eqSPA}
	\end{subequations}
	where $ F: \real^n \to \real $, $ F \in \mathcal{C}^2 $, is strictly convex
	and where $ w: \real^n \to \real^n $ is defined as
	\begin{align}
		w(\dualEq) = - \textstyle{\sum_{ i = 1, i \notin \setEq  }^{n} \dualEq_i \unitVec{i}}
	\end{align}
	with $ \unitVec{i} \in \real^n $ being the $i$th unit vector.
	Let
	\begin{align}
	    \mathcal{M} :=& \big\lbrace (x,\dualEq,\dualIneq) \in \real^{n} \times \real^{n} \times {\realNonNeg^{n}}: \label{eqDefSetOfSaddlePoints}\\
	                 ~& x=x^*, \dualEq_i = 0 \text{ for } i \notin \setEq, \dualIneq_i = 0 \text{ for } i \notin \setIneq, {\textup{and}} \nonumber \\
	                 &\lagrangian(x^*,\dualEq,\dualIneq) \leq \lagrangian(x^*,\dualEq^*,\dualIneq^*) \leq \lagrangian(x,\dualEq^*,\dualIneq^*) \big\rbrace \nonumber 
	\end{align}
	and suppose that \Cref{assMFCQ} holds. 
	Then the set $ \mathcal{M} $ is asymptotically stable
	for \eqref{eqSPA} with region of attraction
	\begin{align}
		\mathcal{R}(\mathcal{M}) \subseteq \big\lbrace (x,\dualEq,\dualIneq) \in \real^{n} \times \real^{n} \times \real^{n}: \dualIneq \in {\realPos^n} \big\rbrace .
		\label{eqDefRegionOfAttraction}
	\end{align}
	\oprocend
\end{lemma}
\begin{remark}
	{Since a point in $ \mathcal{M} $ might 
	as well lie on the boundary of $ \mathcal{R}(\mathcal{M}) $,	
	{one needs to modify the corresponding notions of stability accordingly, by restricting the neighborhoods to the set of admissible initial conditions (cf. \cite{duerr2015thesis}); from now on, we assume that this is understood, without stating it.}} 
	\oprocend
\end{remark}
\begin{remark}
	The {function} $w$ in \eqref{eqSPAb} is usually not included in 
	saddle-point dynamics. Here, it is used to render the
	dynamics of the additional dual variables introduced
	due to the augmentation asymptotically stable.
	It should be noted that the augmentation might lead
	to a significantly larger state vector {for} \eqref{eqSPA}
	compared to the saddle-point dynamics corresponding
	to the original optimization problem \eqref{eqOptimizationProblem}.
	However, it should also be kept in mind that, besides possible
	performance benefits (cf. the discussion after~\Cref{lemmaChoiceSubpath}), the main
	reason for the augmentation is a {significantly} simpler notation and
	it is not crucial for the following methodology to apply (cf. \Cref{remarkAdditionalConstraints}).
	\oprocend
\end{remark}
While \eqref{eqSPA} converges to a solution of \eqref{eqOptimizationProblem},
it is in general not distributed in the aforementioned sense.
Note that if the underlying graph is undirected
and the constraints are only imposed between neighboring agents,
then \eqref{eqSPA} is indeed distributed.
{In the following, we wish} to derive dynamics that ``\emph{approximate}''
those of \eqref{eqSPA} arbitrarily close, {in a sense that will be made precise shortly},
{and are additionally distributed, even when the underlying graph is directed.}
To be more precise, we consider agent dynamics of the form
\begin{subequations}
\begin{align}
    \dot{x}_i^{\seqParam} &= u_{x,i}^{\seqParam}(t,[x_{\mathcal{N}(i)}^{\seqParam},\dualEq_{\mathcal{N}(i)}^{\seqParam},\dualIneq_{\mathcal{N}(i)}^{\seqParam}]) \\
    \dot{\dualEq}_i^{\seqParam} &= u_{\dualEq,i}^{\seqParam}(t,[x_{\mathcal{N}(i)}^{\seqParam},\dualEq_{\mathcal{N}(i)}^{\seqParam},\dualIneq_{\mathcal{N}(i)}^{\seqParam}]) \\
    \dot{\dualIneq}_i^{\seqParam} &= u_{\dualIneq,i}^{\seqParam}(t,[x_{\mathcal{N}(i)}^{\seqParam},\dualEq_{\mathcal{N}(i)}^{\seqParam},\dualIneq_{\mathcal{N}(i)}^{\seqParam}]),
\end{align}\label{eqSPAAgentStructure}
\end{subequations}
where $ i = 1,2,\dots,n$, $ \seqParam \in {\realPos} $ is a parameter and
\begin{align}
    \mathcal{N}(i) := \lbrace j=1,2,\dots,n: g_{ij} \neq 0 \rbrace
\end{align}
is the set of indices of all out-neighboring agents of the $i$th agent. Note that
{the state of the $i$th agent {comprises} of $ ( x_i^\seqParam, \dualEq_i^\seqParam, \dualIneq_i^\seqParam ) $ and}
\eqref{eqSPAAgentStructure} is obviously distributed according to \Cref{defDistributedDynamics}. Our objective
is then to design functions $ u_{x,i}^{\seqParam} $, $ u_{\dualEq,i}^{\seqParam} $,
$ u_{\dualIneq,i}^{\seqParam} $, $ i = 1,2,\dots,n$, parametrized
by $ \seqParam \in {\realPos} $, such that the trajectories
$ \big( x^{\seqParam}(t), \dualEq^{\seqParam}(t), \dualIneq^{\seqParam}(t) \big) $
of \eqref{eqSPAAgentStructure} uniformly converge to the trajectories
$ \big( x(t), \dualEq(t), \dualIneq(t) \big) $ of \eqref{eqSPA} with increasing $ \seqParam $.
To this end, the main idea {of the proposed methodology} is to rewrite the right-hand side of
\eqref{eqSPA} in terms of Lie brackets of \textit{admissible} vector fields,
i.e., vector fields that can be computed locally by the nodes, and 
then employ ideas from geometric control theory to derive
suitable approximations.

\section{Main results}
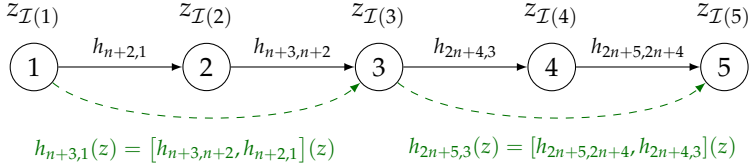
\begin{figure*}[ht]
    \centering
    \begin{minipage}{0.6\textwidth}\vspace{0pt}
	\begin{tikzpicture}[>=latex]
		\node[circle,draw,name=node1] at (0,0) {$1$};
		\node[name=labelNode1,above=0.2em of node1] {$ \state_{\mathcal{I}(1)} $};
		\node[circle,draw,name=node2,right=1.65cm of node1] {$2$};
		\node[name=labelNode2,above=0.2em of node2] {$ \state_{\mathcal{I}(2)} $};
		\node[circle,draw,name=node3,right=1.65cm of node2] {$3$};
		\node[name=labelNode3,above=0.2em of node3] {$ \state_{\mathcal{I}(3)} $};
		\node[circle,draw,name=node4,right=1.65cm of node3] {$4$};
		\node[name=labelNode4,above=0.2em of node4] {$ \state_{\mathcal{I}(4)} $};
		\node[circle,draw,name=node5,right=1.65cm of node4] {$5$};
		\node[name=labelNode5,above=0.2em of node5] {$ \state_{\mathcal{I}(5)} $};
		\draw[->] (node1) -- node[scale=0.8,anchor=south] {$ h_{n+2,1} $} (node2);
		\draw[->] (node2) -- node[scale=0.8,anchor=south] {$ h_{n+3,n+2} $} (node3);
		\draw[->] (node3) -- node[scale=0.8,anchor=south] {$ h_{2n+4,3} $} (node4);
		\draw[->] (node4) -- node[scale=0.8,anchor=south] {$ h_{2n+5,2n+4} $} (node5);
		\draw[->,dashed,color=fictitiousColor] (node1) [out=-40, in=-140,looseness=0.6] to (node3); 
		\node[below=0.5cm of node2,scale=0.8,color=fictitiousColor,xshift=-10pt] { $ h_{n+3,1}(\state) = \big[ h_{n+3,n+2},h_{n+2,1} \big](\state) $ };
		\draw[->,dashed,color=fictitiousColor] (node3) [out=-40, in=-140,looseness=0.6] to (node5);
		\node[below=0.5cm of node4,scale=0.8,color=fictitiousColor,xshift=10pt] { $ h_{2n+5,3}(\state) = [ h_{2n+5,2n+4}, h_{2n+4,3} ](\state) $ };
	\end{tikzpicture}
    \end{minipage}
    \begin{minipage}{0.39\textwidth}\vspace{0pt}
\caption{A communication structure with $n=5$ nodes is depicted. The arrows
indicate to which agent state some agent has access to, e.g., agent $1$ has access to the state
of agent $2$ {given by $ \state_{\mathcal{I}(2)} = [ x_2, \dualEq_2, \dualIneq_2 ]^\top $}
but not the other way round.
The dotted green arrow shows a fictitious edge
with associated vector fields created by Lie brackets of admissible 
vector fields.}
\label{figGraphChain}
\end{minipage}
\end{figure*}
Consider the saddle-point dynamics \eqref{eqSPA}.
As a first step, we separate the right-hand side
into admissible and non-admissible vector fields,
where admissible refers to the part of the 
dynamics that can be computed locally by the nodes.
For the ease of presentation{,} we assume in the following
that the constraints of agent $i$ are only imposed to
its out-neighboring agents, i.e., we impose
the following assumption on the constraints:
\begin{assumption}\label{assConstraints2}
For $ \eqConstrLAgent_i = [ \eqConstrLAgent_{i1}, \dots, \eqConstrLAgent_{in} ] $,
	$ \ineqConstrLAgent_i = [ \ineqConstrLAgent_{i1}, \dots, \ineqConstrLAgent_{in} ] $,
	$ i = 1,2,\dots,n $, we have for each $j=1,\dots,n$,
	that 	$ \eqConstrLAgent_{ij} \neq 0 $ or $ \ineqConstrLAgent_{ij} \neq 0 $
	only if $ \elLaplacian_{ij} \neq 0 $. \oprocend
\end{assumption}
In other words, we thereby assume that the constraints
match the communication topology induced by the graph\footnote{It should be noted that the following results can be
extended to problems where this assumption does not
hold, cf.{~\cite{mic2018extensions},} \Cref{remarkConstraintsNotBetweenNeighbors}
as well as the example in~\Cref{secSimExample}.}.
Under this assumption, the right-hand side of \eqref{eqSPAb}, \eqref{eqSPAc}
is admissible, while parts of the right-hand side of
\eqref{eqSPAa} are not. {Note that} the gradient
of $F$ is admissible, since $F$ is a separable function;
{the remaining terms, however, are not necessarily
admissible, since the underlying communication graph
is directed.}
Now, for $ {\eqConstrL} = [ \eqConstrLAgent_{ij} ] $, $ {\ineqConstrL} = [ \ineqConstrLAgent_{ij} ] $,
we define the admissible part of $ \eqConstrL^\top, \ineqConstrL^\top $ as
\begin{align}
    \tilde{\eqConstrL}_{\textup{adm}} 
    &=
    \sum\limits_{i=1}^{n} \sum\limits_{j=1}^{n} \sign( \vert \elLaplacian_{ij} \vert ) 
    \eqConstrLAgent_{ji} \unitVec{i} \unitVec{j}^\top, \\
    \tilde{\ineqConstrL}_{\textup{adm}} 
    &=
    \sum\limits_{i=1}^{n} \sum\limits_{j=1}^{n} \sign( \vert \elLaplacian_{ij} \vert )  
    \ineqConstrLAgent_{ji} \unitVec{i} \unitVec{j}^\top,
\end{align}
where {$ \sign: \real \to \lbrace -1,0,1 \rbrace $} is the sign function and $ \unitVec{i} $ is the $i$th unit vector.
Observe that $ \tilde{\eqConstrL}_{\textup{adm}} $, $ \tilde{\ineqConstrL}_{\textup{adm}} $ 
correspond to the admissible part of $ \eqConstrL^\top $ and $ \ineqConstrL^\top $,
respectively. We then let
\begin{align}
     \tilde{\eqConstrL}_{\textup{rest}} &= \eqConstrL^\top - \tilde{\eqConstrL}_{\textup{adm}},
     \quad
     \tilde{\ineqConstrL}_{\textup{rest}} = \ineqConstrL^\top - \tilde{\ineqConstrL}_{\textup{adm}},
\end{align}
{and define the state of~\eqref{eqSPA} as}
\begin{align}
	\state := [ x^\top, \dualEq^\top, \dualIneq^\top ]^\top \in \real^{3n}.
	\label{eqDefCompleteState}
\end{align}
Hence, we can write the saddle-point dynamics \eqref{eqSPA}
as
\begin{align}
	{\dot{\state}}
	&=
	{f_{\textup{adm}}(\state)}
	+
	\begin{bmatrix}
	- \tilde{\eqConstrL}_{\textup{rest}} \dualEq - \tilde{\ineqConstrL}_{\textup{rest}} \dualIneq \\
	0 \\
	0
	\end{bmatrix}, \label{eqSPASeparated}
\end{align}
where $ f_{\textup{adm}}:\real^{3n} \to \real^{3n} $ is defined as
\begin{align}
	{f_{\textup{adm}}(\state)} 
	=
	\begin{bmatrix} 
	-\nabla F(x) - \tilde{\eqConstrL}_{\textup{adm}} \dualEq - \tilde{\ineqConstrL}_{\textup{adm}} \dualIneq \\
	\eqConstrL x - \eqConstrR + w(\dualEq) \\
	\textup{diag}(\dualIneq) \big( \ineqConstrL x - \ineqConstrR \big)
	\end{bmatrix}.
	\label{eqDefAdmissiblePart}
\end{align}
{Here}, $ f_{\textup{adm}} $ is admissible whereas
the second {term on the right-hand side of}~\eqref{eqSPASeparated} is not. The essential idea to derive
suitable distributed approximations is to rewrite
the non-admissible part in terms of Lie brackets of
admissible vector fields{; we will elaborate on
this in what follows next}.

\subsection{Rewriting the non-admissible vector fields}
{We first define the index set}
\begin{align}
	\mathcal{I}(i) := \lbrace i, n+i, 2n+i \rbrace,
	\label{eqDefIndexSetAgent}
\end{align}
{where} $ i = 1,2,\dots,n $, associating the components of $ \state $ to the
$i$th agent, i.e., {$ \state_{\mathcal{I}(i)} $ is the state of agent~$i$.}
We then define a set of vector fields $ h_{i,j}: \real^{3n} \to \real^{3n} $,
$ i,j = 1,2,\dots,3n $, as
\begin{align}
	h_{i,j}(\state) = \state_i \unitVec{j},  \label{eqDefAdmissibleVectorFields}
\end{align}
where $ \unitVec{j} \in \real^{3n} $ is the $j$th unit vector.
Observe that $ h_{i,j} $ is an admissible vector field if and only 
if there exist $ \ell, k $ such that $ i \in \mathcal{I}(\ell) $,
$ j \in \mathcal{I}(k) $ and $ \elLaplacian_{k\ell} \neq 0 $.
Before we present a general construction rule, let us first illustrate
the main idea by means of a simple example.
\begin{example}\label{exampleRewritingBrackets}
Consider the graph shown in~\Cref{figGraphChain} with $n=5$ nodes.
	Let $h_{i,j}$ be defined as in \eqref{eqDefAdmissibleVectorFields}
	and observe that $ h_{n+3,n+2}, h_{n+2,1} $ are admissible.
	Consider the Lie bracket
	\begin{align}
		&\big[ h_{n+3,n+2}, h_{n+2,1} \big](\state) \nonumber \\
		=~&
		\unitVec{1} \unitVec{n+2}^\top \state_{n+3} \unitVec{n+2} - \unitVec{n+2} \unitVec{n+3}^\top \state_{n+2} \unitVec{1} \nonumber \\
		=~& 
		\state_{n+3} \unitVec{1},
	\end{align}
	which, according to \eqref{eqDefAdmissibleVectorFields}, is equal
	to $ h_{n+3,1}(\state) $, i.e., a non-admissible vector field. Given the graphical representation
	in \Cref{figGraphChain}, this can be interpreted as a ``fictitious'' 
	edge from agent $1$ to agent $3$, generated by the Lie bracket of two 
	admissible vector fields. This observation is of key importance
	in the rest of the paper. More generally, we can observe that
	\ifthenelse{\boolean{longVersion}}{
	\begin{align}
		\big[ h_{i,j}, h_{j,k} \big](\state) = h_{i,k}(\state),
	\end{align}}{
	{$\big[ h_{i,j}, h_{j,k} \big](\state) = h_{i,k}(\state)$, }}
	for any $ i,j,k = 1,2,\dots,3n $.
	\oprocend
\end{example}
Next, we generalize this idea. Let $ p_{ij} = \pathLeft v_{i_1} \pathSep \ldots \pathSep v_{i_{r}} \pathRight $
be a path in $ \mathcal{\laplacian} = (\mathcal{V},\mathcal{E}) $ 
from node $v_i$ to node $v_j$, i.e. $ i=i_1 $, $ j=i_{r} $,
$ v_{i_1}, \ldots, v_{i_r} \in \mathcal{V} $, $ r \geq 2 $,
and let $ \length(p_{ij}) = r-1 $ denote its length.
We now, recursively, define a mapping $ \recBracket{k_1,k_2}{} $, $ k_1,k_2 = 1,2,\dots,3n$, 
from a given path $ p_{ij} $ in $ \mathcal{G} $ 
to the set of vector fields on $ \real^{3n} $:
\begin{itemize}[leftmargin=*]
	\item for $ \ell(p_{ij}) = 1$, we define 
	\begin{align}
		\recBracket{k_1,k_2}{}(p_{ij}) = h_{k_1,k_2}.
		\label{eqDefRecursionBStart}
	\end{align}
	\item {for $ \ell(p_{ij}) \geq 2 $, we define
	\begin{align}
		\recBracket{k_1,k_2}{}(p_{ij}) = [ \recBracket{k_1,s}{}(q^c), \recBracket{s,k_2}{}(q) ],
		\label{eqDefRecursionBStep}
	\end{align}
	where $ q $ is any subpath in $ \subpath_{i \bullet}(p_{ij}) $ and 
	$ s \in \mathcal{I}(\tail(q)) $.}
\end{itemize}
Observe that $ \recBracket{k_1,k_2}{} $ is independent of the path $ p_{ij} $ according
	to the definition~\eqref{eqDefRecursionBStart}. However, the path comes into
	play when it gets to choosing $k_1,k_2$ such that the resulting Lie bracket 
	is a Lie bracket of admissible vector fields, cf.~\Cref{lemmaConstructionBrackets}.
Using~\eqref{eqDefRecursionBStart}, \eqref{eqDefRecursionBStep}, we next state a result that
extends the ideas from \Cref{exampleRewritingBrackets}{; a proof is provided in~\Cref{secAppendixProofConstructionBrackets}.}
\begin{lemma}\label{lemmaConstructionBrackets}
Consider a directed graph $ \mathcal{G} = ( \mathcal{V}, \mathcal{E} ) $ of $n$ nodes.  
	Let $p_{ij}$ be a path between $v_i$ and $v_j$, $ v_i,v_j \in \mathcal{V}$, 
	{and let $ \recBracket{k_1,k_2}{} $ be defined as in~\eqref{eqDefRecursionBStart},~\eqref{eqDefRecursionBStep}.}
	Then, if \mbox{$k_1 \neq k_2$}, we have for all $ \state \in \real^{3n} $
	\begin{align}
		\recBracket{k_1,k_2}{}\big(p_{ij}\big)(\state) = \state_{k_1} \unitVec{k_2} = h_{k_1,k_2}(\state),
		\label{eqLemmaBracket}
	\end{align}
	and, if $ k_1 \in \mathcal{I}\big( \tail(p_{ij}) \big) $, $ k_2 \in \mathcal{I}\big( \head(p_{ij}) \big) $, 
	then	$ \recBracket{k_1,k_2}{}(p_{ij}) $ is a Lie bracket of admissible vector fields.
	\oprocend
\end{lemma}
\begin{remark}\label{remarkAdditionalConstraints}
The same result holds true if we drop the assumption that
	each agent has exactly one equality and one inequality
	constraint, since this only leads to a {reformulation} of the index sets
	$ \mathcal{I}(i) $, $ i = 1,2,\dots,n $. Interestingly,
	additional constraints also introduce additional degrees
	of freedom in rewriting the non-admissible vector fields{,}
	since the index set $ \mathcal{I}(\tail(q)) $ grows.
	\oprocend
\end{remark}
\begin{remark}
{It is worth pointing out that} {admissible} vector fields of the form \eqref{eqDefAdmissibleVectorFields}
	are not the only ones that can be used to rewrite {(linear)} non-admissible
	vector fields in terms of Lie brackets of admissible vector fields.
	In fact, {as discussed in~\cite{mic2018extensions} in detail,}
	there exists a whole class of {admissible} vector fields which can be
	employed for this purpose. {Similar as in~\cite{grushkovskaya2018class},
	a different choice can positively affect the approximation quality of the 
	resulting distributed algorithm.}
	\oprocend
\end{remark}
While \Cref{lemmaConstructionBrackets} holds for any directed
path in $ \mathcal{G} $, {from now on} we use the
shortest path {as it} leads to iterated Lie brackets 
of smallest degree. We do not discuss how to compute
the paths here since this is a problem on its own but refer
the reader to standard algorithms, see, e.g., \cite{cormen2009introduction}.
Further, the choice of subpath
and the state index $s$ in the recursion~\eqref{eqDefRecursionBStep} 
is arbitrary as well. In~\Cref{lemmaChoiceSubpath} in~\Cref{secInputConstruction}, we provide a particular choice
that turns out to be beneficial in the construction
of the approximating input sequences. 
The next result is an immediate consequence of
\Cref{lemmaConstructionBrackets}.
\begin{proposition}\label{propRewriting}
Suppose that \Cref{assConstraints2} holds {and that}
	$ \mathcal{G} = ( \mathcal{V}, \mathcal{E} ) $ is strongly connected.  
	For all $ i,j = 1,\dots,n $, let $ p_{ij} $, denote 
	a path from node $ v_i $ to node $ v_j $, {where} $ v_i, v_j \in \mathcal{V} $.
	Then, with $ \state = [ x^\top, \dualEq^\top, \dualIneq^\top ]^\top $, the dynamics \eqref{eqSPASeparated} can
	equivalently be written as
	\begin{align}
	\begin{split}
		\dot{\state}
		&=
		f_{\textup{adm}}(\state)
		-\sum\limits_{i=1}^n \sum\limits_{j=1}^n \tilde{\MakeLowercase{\eqConstrL}}_{\textup{rest},ij} \, \recBracket{n+j,i}{}\big( p_{ij} \big) ( \state ) \\
		&-\sum\limits_{i=1}^n \sum\limits_{j=1}^n \tilde{\MakeLowercase{\ineqConstrL}}_{\textup{rest},ij} \recBracket{2n+j,i}{}\big( p_{ij} \big) ( \state ) 
		\end{split}\label{eqSPArewritten}
	\end{align}
	and the right-hand side is a linear combination of Lie brackets
	of admissible vector fields.
	\oprocend
\end{proposition}
\begin{remark}\label{remarkConstraintsNotBetweenNeighbors}
If \Cref{assConstraints2} does not hold the terms $ [0, f_{\textup{adm},2}(\state), 0 ]^\top $, 
	$ [0, 0, f_{\textup{adm},3}(\state) ]^\top $
	may no longer be admissible. While $ [0, f_{\textup{adm},2}(\state), 0 ]^\top $
	can be rewritten using \Cref{lemmaConstructionBrackets}, for
	$ [0, 0, f_{\textup{adm},3}(\state) ]^\top $
	different construction techniques are required, since $ f_{\textup{adm},3} $
	is bilinear {as a function of $ x $ and $ \dualIneq $}. However, it should be noted that
	it is {still} possible to rewrite these terms {by means of admissible vector fields, see~\cite{mic2018extensions}.}
	\oprocend
\end{remark}
\begin{remark}
In general, {having} a strongly connected graph is sufficient but not necessary. 
	In fact, it is sufficient that there exists a path from node $v_i$
	to node $v_j$ for all $ i,j $ such that $ \tilde{\MakeLowercase{\eqConstrL}}_{\textup{rest},ij} \neq 0 $ 
	or $ \tilde{\MakeLowercase{\ineqConstrL}}_{\textup{rest},ij} \neq 0 $.
	\oprocend
\end{remark}
Now that we have rewritten the non-admissible
vector fields in terms of iterated Lie brackets of
admissible vector fields, there is still the issue 
of {generating} suitable functions $ u_{x,i}^{\seqParam},
u_{\dualEq,i}^{\seqParam}, u_{\dualIneq,i}^{\seqParam} $
to be addressed. We will {study} this in
the next section and provide a result
on how \eqref{eqSPAAgentStructure}
and \eqref{eqSPArewritten} are related in terms of
their stability properties under a suitable
choice of the input functions.

\subsection{Construction of distributed control laws}\label{secInputConstruction}
{Our main objective in this section is to} elaborate on how to construct suitable input
functions $ u_{x,i}^{\seqParam}, u_{\dualEq,i}^{\seqParam}, u_{\dualIneq,i}^{\seqParam} $
such that the trajectories of \eqref{eqSPAAgentStructure} uniformly
converge to those of \eqref{eqSPArewritten} {as we increase}
$ \seqParam $. The following procedure is based on the
results presented in \cite{liu1997approximation}, \cite{sussmann1991limits}, \cite{liu1997averaging}.
In \cite{liu1997averaging}, the relation between the
trajectories of a system of the form
\begin{align}
	\dot{\state}^{\seqParam} = f_0(\state^{\seqParam}) + \sum\limits_{k=1}^{M} \phi_k(\state^{\seqParam}) U_k^{\seqParam}(t), \quad \state^{\seqParam}(0) = \state_0,
	\label{eqOriginalSystemLiu}
\end{align}
{where} $ f_0, \phi_k: \real^N \to \real^N $, $ U_k^{\seqParam}: \real \to \real $,
$ \state_0 \in \real^N $ and the trajectories of an associated \emph{extended system} 
\begin{align}
	\dot{\state} = f_0(\state) + \sum\limits_{B \in \mathcal{B}} v_B B(\state), \quad \state(0) = \state_0,
	\label{eqExtendedSystemLiu}
\end{align}
is studied, where $ \mathcal{B} $ is a {finite} set of Lie brackets
of the vector fields $ \phi_k $, $ k = 1,\dots, M $,
and $ v_B \in \real $ is the corresponding coefficient.
In our setup, \eqref{eqSPAAgentStructure} will play the role of
\eqref{eqOriginalSystemLiu} with $ \phi_k $ being the admissible vector
fields and \eqref{eqSPArewritten} plays the role of 
\eqref{eqExtendedSystemLiu} with $ \mathcal{B} $ being
the set of Lie brackets of admissible vector fields required
to rewrite the non-admissible vector fields.
It is shown in \cite{liu1997averaging} that, under a suitable choice of the input
functions $ U_k^{\seqParam} $, the solutions of \eqref{eqOriginalSystemLiu}
uniformly converge to those of \eqref{eqExtendedSystemLiu}
on compact time intervals for increasing $ \seqParam $, i.e.,
for each {$z_0 \in \real^N$}, for each $ \varepsilon > 0 $ and for each $T \geq 0 $, there exists
$ \seqParam^* > 0 $ such that for all $ \seqParam > \seqParam^* $ {and $ t \in [0,T] $ we have that}
\ifthenelse{\boolean{longVersion}}{
\begin{align}
	\Vert \state(t) - \state^{\seqParam}(t) \Vert \leq \varepsilon.
\end{align}
}
{$ \Vert \state(t) - \state^{\seqParam}(t) \Vert \leq \varepsilon$.}
An algorithm {for constructing} suitable input functions $ U_k^{\seqParam} $ 
that fulfill these assumptions is presented in~\cite{liu1997approximation} {as well as in a brief version in~\cite{sussmann1991limits}}; we will follow
this idea in here, however, given that in~\cite{liu1997approximation}
the input functions are not given in explicit form, we exploit
the special structure of the admissible vector fields in order
{to simplify this procedure and arrive at} explicit formulas for
{a large class of scenarios applicable to our work.}

\subsubsection{Writing the Lie brackets in terms of a P.~Hall basis}
The algorithm presented in \cite{liu1997approximation} requires
the brackets {used} in~\eqref{eqExtendedSystemLiu} to be brackets
in a {so-called P.~Hall basis}; we need to {``project''} the
brackets in~\eqref{eqSPArewritten} to such a basis, {in the sense that will be made precise shortly}.
We first recall the definition of a P. Hall basis;
{we let $ \degree(B) $ denote the degree of a bracket $B$.}
\begin{definition}[P.~Hall basis of a Lie algebra]\label{defPHallBasis}
Let $ \Phi = \lbrace \phi_1, \phi_2, \dots, \phi_{{M}} \rbrace $ be a set of {smooth} vector fields.
	A P.~Hall basis $\mathcal{PH}(\Phi) = ( \mathbb{P},\prec ) $
	of the Lie algebra {generated by} $ \Phi $ is a set $ \mathbb{P} $ of 
	brackets equipped with a total ordering $ \prec $ that fulfills the
	following properties: \\[-2em]
	\begin{enumerate}[leftmargin=*,label={[\textup{PH}\arabic*]}]
		\item Every $\phi_k${, $k=1,2,\dots,M$,} is in $\mathbb{P}$. \label{itemPH1}
		\item $ \phi_k \prec \phi_j $ if and only if $ k < j $. \label{itemPH2}
		\item If $ B_1, B_2 \in \mathbb{P} $ and $ \degree(B_1) < \degree(B_2) $, then $ B_1 \prec B_2 $. \label{itemPH3}
		\item Each $ B = [B_1,B_2] \in \mathbb{P} $ if and only if \label{itemPH4}
		\begin{enumerate}[label={{[PH4.\alph*]}}]
		 \item $ B_1, B_2 \in \mathbb{P} $ and $ B_1 \prec B_2 $ \label{itemPH4a}
		 \item either $ \degree(B_2) = 1 $ or $ B_2 = [B_3,B_4] $  for some $ B_3, B_4 $ such that $ B_3 \preceq B_1 $. \label{itemPH4b} \oprocend
		\end{enumerate}
	\end{enumerate}
\end{definition}
{
\begin{remark}
It is understood that a P.~Hall basis is well-defined only for 
	formal brackets of indeterminates but not for Lie brackets of vector fields.
	In particular, in~\ref{itemPH3} and~\ref{itemPH4}, for Lie brackets
	the degree as well as the left and right factors $ B_1 $ and $ B_2 $ are not uniquely
	defined, see also~\Cref{secPreliminaries}. 
	For the purpose of a clearer presentation we avoid this formal overhead
	accepting this 	abuse of notation and assume that $B$ is
	interpreted as a formal bracket in~\ref{itemPH3}, \ref{itemPH4}.
	The interested reader is referred to
	\ifthenelse{\boolean{longVersion}}{
	\Cref{secAppendixFormalBrackets}
	}
	{the extended version~\cite{mic2017distributedArxiv}}
	for some more details on this subject.
\end{remark}
}
Note that
\ref{itemPH2} is usually not included in the definition 
of a P.~Hall basis, but it is common to include it for the 
approximation problem at hand. 
Moreover, the construction rule \ref{itemPH4} 
ensures that no brackets are included in the basis that are
related to other brackets in the basis by the Jacobi identity
{or skew-symmetry; thus} the brackets are in this sense independent.
However, this does not mean that, when evaluating the brackets, the resulting vector fields are
independent, which we will exploit later.
It is as well worth mentioning that the ordering fulfilling the
properties \ref{itemPH1} - \ref{itemPH4} is in general
not unique, i.e., for a given set of vector fields
$ \Phi $, there {may} exist several P.~Hall bases.

{Let us now return to our setup. Let $ \Phi $ be given by the set}
of admissible vector fields defined as
\begin{align}	
	\Phi := \big\lbrace h_{i,j} :~& \exists k_1,k_2 { \in \lbrace 1,2,\dots,n \rbrace} \text{ such that } i \in \mathcal{I}(k_1), \nonumber \\
	&j \in \mathcal{I}(k_2),  \elLaplacian_{k_2k_1} \neq 0 \big\rbrace,
	\label{eqDefSetOfAdmissibleVectorFields}
\end{align}
where $ h_{i,j} $ is defined in \eqref{eqDefAdmissibleVectorFields}.
Every bracket in the set of Lie brackets of admissible vector fields $ \mathcal{B} $ can then be 
{projected onto some P. Hall basis $ \mathcal{PH}(\Phi) $, i.e., be uniquely} 
written as a linear combination of elements of {$ \mathcal{PH}(\Phi) $} 
by successively resorting the 
brackets, making use of skew-symmetry and the Jacobi identity{, cf.~\Cref{remarkProjectionPHall} for an example}.
{Such a projection
algorithm is for example given in \cite{reutenauer2003free}
and in the following we let for any $ B \in \lieBr(\Phi) $
\begin{align}
	\resort_{\mathbb{P}}(B) = \sum\limits_{\tilde{B} \in \mathbb{P}} \theta_{\tilde{B}} \tilde{B}
	\label{eqDefProjectionPHall}
\end{align}
denote the unique representation of $B$ in terms of brackets
from a P.~Hall basis $ \mathcal{PH}(\Phi) = (\mathbb{P},\prec) $.}
However, for brackets of higher degree, {finding this representation} 
might be tedious and result{s}
in a large number of brackets {$\tilde{B}$}; {we hence propose an alternative approach.
Instead of resorting the complete
brackets appearing in \eqref{eqSPArewritten},} we suggest
to reduce the resorting {to brackets of low degree} by a proper choice
of the subpaths in the construction procedure
presented in \Cref{lemmaConstructionBrackets}.
The main idea is to choose the subpath $q$
in \eqref{eqDefRecursionBStep} in such a way that,
in each recursion step, the degree of the left factor
of the bracket is strictly smaller than the degree
of the right factor and such that the degree of the
left factor of the right factor is smaller than that
of the left factor of the original bracket such that
\ref{itemPH4a} and \ref{itemPH4b} are automatically
fulfilled. Since the degree directly corresponds to
the length of the subpath this can be achieved by 
choosing the subpath appropriately\ifthenelse{\boolean{longVersion}}{, see also \Cref{figIllustrationSubpaths}}{{, see also~\cite{mic2017distributedArxiv} for an illustration}}. We make this idea more precise
in the following Lemma.
\begin{lemma}\label{lemmaChoiceSubpath}
Consider a directed graph $ \mathcal{G} = ( \mathcal{V}, \mathcal{E} ) $ of $n$ nodes.  
	Let the set of admissible vector fields be defined according to \eqref{eqDefSetOfAdmissibleVectorFields}.
	Let some P.~Hall basis \mbox{$ \mathcal{PH}(\Phi) = (\mathbb{P},\prec) $} be given
	and let {$\resort_{\mathbb{P}}(B)$ denote the unique representation of $B$ 
	in terms of brackets in $\mathbb{P} $, cf.~\eqref{eqDefProjectionPHall}.}
	Let $ p_{i_1 i_r} $ be a path from node
	$ v_{i_1} \in \mathcal{V} $ to node $ v_{i_r} \in \mathcal{V} $
	and define
	\begin{align}
		&\recBracketPHall{k_1,k_2}{}(p_{i_1 i_r}) \label{eqRecPHall} \\
		&=
		\begin{cases}
			\recBracket{k_1,k_2}{}(p_{i_1 i_r})                                    & \text{if } \length(p_{i_1 i_r}) = 1, \\
			\resort_{{\mathbb{P}}}\big( [ \recBracket{k_1,s}{}( q^c ), \recBracket{s,k_2}{}( q )] \big) \hspace*{-0.6em}  & \text{if } \length(p_{i_1 i_r}) = 2,3,4,6, \\
			[ \recBracketPHall{k_1,s}{}( q^c ), \recBracketPHall{s,k_2}{}( q )]    & \text{{otherwise,}}
		\end{cases} \nonumber
	\end{align}
	where 
	\begin{align}
		s &= 
		\begin{cases}
			n + i_{\theta(p_{i_1 i_r})} & \text{if } 1 \leq k_1 \leq 2n \\
			{2n} + i_{\theta(p_{i_1 i_r})} & \text{if } 2n+1 \leq k_1 \leq 3n \\
		\end{cases} \label{eqChoiceS} \\
		q &= p_{i_1  i_{\theta(p_{i_1 i_r})} } \in \subpath_{i_1 \bullet} (p_{i_1 i_r})
		\label{eqChoiceSubpath} \\
		\theta(p_{i_1 i_r}) &= 
		\begin{cases}
			\tfrac{1}{2}\length(p_{i_1 i_r}) + 1                 & \text{if }  \length(p_{i_1 i_r}) = 2,4, \\
			\lfloor \tfrac{1}{2}\length(p_{i_1 i_r}) \rfloor + 2 & \text{otherwise},
		\end{cases}
		\label{eqChoiceTheta}
	\end{align}
	with $ \lfloor a \rfloor $ being the largest integer value {less or equal} than $ a \in {\realNonNeg} $. 
	Then $ \recBracketPHall{k_1,k_2}{}(p_{i_1i_r})(\state) = \recBracket{k_1,k_2}{}(p_{i_1i_r})(\state) $
	for all $ \state \in \real^{3n} $ and $ \recBracketPHall{k_1,k_2}{}(p_{i_1i_r}) \in \mathbb{P} $
	for all $ k_1 \in \mathcal{I}(\tail(p_{i_1 i_r})),  k_2 \in \mathcal{I}(\head(p_{i_1 i_r})) $.
	\oprocend
\end{lemma}
{A proof is given in \Cref{secAppendixProofChoiceSubpath}. Equation~\eqref{eqRecPHall} and the choice of $s, q$ from~\eqref{eqChoiceS}, \eqref{eqChoiceSubpath}
can be interpreted as follows: A bracket corresponding to a path $p_{i_1 i_r}$
of length larger than one is generated by dividing the path into two complementing subpaths
$q$ and $q^c$, where~\eqref{eqChoiceTheta} ensures that the resulting brackets
have the desired properties~\ref{itemPH4}. The cases where these properties are not ensured by that choice, i.e.,
$ \length(p_{i_1 i_r}) \in \lbrace 2,3,4,6 \rbrace$, are handled separately. 
Further, $s$ corresponds, roughly speaking, to the element of the complete state
vector over which the information is passed. As it turns out in the design
of the approximating inputs, this also corresponds to the components of the complete
state in which the perturbing inputs are injected. It is worth pointing out, as become clear in the proof, 
that the aforementioned result is independent of the choice of $s$ as given in~\eqref{eqChoiceS}; in fact, any
{$ s \in \mathcal{I}\big(i_{\theta(p_{i_1 {i_r}})}\big) = \lbrace i_{\theta(p_{i_1 i_r})}, n+ i_{\theta(p_{i_1 i_r})}, 2n + i_{\theta(p_{i_1 i_r})} \rbrace $}
can be taken. The specific choice \eqref{eqChoiceS} has advantages that will be made clear later.
Observe that the degrees of freedom for $s$ increase with the number 
of constraints of each agent. In particular, it might as well
happen that there is no degree of freedom if we do not augment
the optimization problem \eqref{eqOptimizationProblem}. 
}
\begin{remark}\label{remarkProjectionPHall}
It should be noted that the projection can be computed easily in the 
	given case. To this end, first notice that -- by the choice of subpaths --
	for $ \length(p_{i_1 i_r}) = 2,3 $, the brackets admit the following structure
	\begin{align}
		\recBracket{k_1,k_2}{}(p_{i_1 i_r})
		=
		\begin{cases}
			[ \phi_{a_1}, \phi_{a_2} ] & \text{if } \length(p_{i_1 i_r}) = 2 \\
			\big[ \phi_{a_1}, [ \phi_{a_2}, \phi_{a_3} ] \big] & \text{if } \length(p_{i_1 i_r}) = 3 
		\end{cases}
	\end{align}
	for some $ a_{1/2/3} \in {\natPos} $ depending on $k_1,k_2,p_{i_1 i_r} $,
	where $ \phi_{a_i} \in \Phi $, $ i = 1,2,3 $. For {such} brackets,
	the projection on the P.~Hall basis $ \mathcal{PH}{(\Phi) = (\mathbb{P},\prec)} $ is easily computed 
	making use of skew-symmetry and the Jacobi-identity and we obtain
	\begin{align}
		\resort_{{\mathbb{P}}} \big( [ \phi_{a_1}, \phi_{a_2} ] \big) 
		&=
		\begin{cases}
			\hphantom{-} [ \phi_{a_1}, \phi_{a_2} ] & \text{if } a_1 < a_2, \\
			- [ \phi_{a_2}, \phi_{a_1} ] & \text{if } a_1 > a_2,
		\end{cases} \label{eqProjectionPHallDegree2}
	\end{align}
	and
	\begin{align}
		&\resort_{{\mathbb{P}}} \big( \big[ \phi_{a_1}, [ \phi_{a_2}, \phi_{a_3} ] \big] \big) = \label{eqProjectionPHallDegree3} \\
		&\begin{cases}
			\big[ \phi_{a_2}, [ \phi_{a_1}, \phi_{a_3} ] \big] - \big[ \phi_{a_3}, [ \phi_{a_1}, \phi_{a_2} ] \big] & \text{if } a_1 = \min\limits_{i = 1,2,3} a_i, \\
			\big[ \phi_{a_1}, [ \phi_{a_2}, \phi_{a_3} ] \big] & \text{if } a_2 = \min\limits_{i = 1,2,3} a_i, \\
			- \big[ \phi_{a_1}, [ \phi_{a_3}, \phi_{a_2} ] \big] & \text{if } a_3 = \min\limits_{i = 1,2,3} a_i.
		\end{cases} \nonumber
	\end{align}
	{Note that the brackets have been resorted in such a way that 
	the brackets on the right hand side of~\eqref{eqProjectionPHallDegree2},
	\eqref{eqProjectionPHallDegree3} fulfill~\ref{itemPH3},~\ref{itemPH4}
	when interpreted as formal brackets.}
	In the same manner, for $ \length(p_{i_1 i_r}) = 4,6 $, we have
	\begin{flalign}
		\recBracket{k_1,k_2}{}(p_{i_1 i_r})
		&=
		\begin{cases}
			[ B_{a_1}, B_{a_2} ] & \text{if } \length(p_{i_1 i_r}) = 4, \\
			\big[ B_{a_1}, [ B_{a_2}, B_{a_3} ] \big] & \text{if } \length(p_{i_1 i_r}) = 6,
		\end{cases} \hspace*{-2em} & \label{eqProjectionPHallDegree4_6}
	\end{flalign}
	where the $ B_{a_i} $ are Lie brackets of the $ \phi_i $ with $ \degree(B_{a_i}) = 2 $, $ i = 1,2,3 $.
	The projection is then done by first projection the inner brackets $ B_{a_i} $ 
	on the P.~Hall basis using \eqref{eqProjectionPHallDegree2} and then
	resorting $ \recBracket{k_1,k_2}{}(p_{i_1 i_r}) $ as in \eqref{eqProjectionPHallDegree2}, \eqref{eqProjectionPHallDegree3}.
	\oprocend
\end{remark}

\ifthenelse{\boolean{longVersion}}{
\begin{figure*}
	\begin{center}
		\begin{tikzpicture}[>=latex]
			\node[name=nodePathLabel] {$p_{1 6}$};
			\node[name=node6,right=1cm of nodePathLabel,draw,circle] {$6$};
			\node[name=node5,right=1cm of node6,draw,circle] {$5$};
			\node[name=node4,right=1cm of node5,draw,circle] {$4$};
			\node[name=node3,right=1cm of node4,draw,circle] {$3$};
			\node[name=node2,right=1cm of node3,draw,circle] {$2$};
			\node[name=node1,right=1cm of node2,draw,circle] {$1$};
			\draw[->] (node1) -- (node2);
			\draw[->] (node2) -- (node3);
			\draw[->] (node3) -- (node4);
			\draw[->] (node4) -- (node5);
			\draw[->] (node5) -- (node6);
			\node[name=nodeSubPath1Label,below=1cm of nodePathLabel.east,anchor=east] {$q = p_{14} \in \subpath_{1 \bullet}(p_{16})$};
			\node[name=node1SubPath1,draw,circle] at (node1 |- nodeSubPath1Label) {$1$};
			\node[name=node2SubPath1,left=1cm of node1SubPath1,draw,circle] {$2$};
			\node[name=node3SubPath1,left=1cm of node2SubPath1,draw,circle] {$3$};
			\node[name=node4SubPath1,left=1cm of node3SubPath1,draw,circle] {$4$};
			\draw[->] (node1SubPath1) -- (node2SubPath1);
			\draw[->] (node2SubPath1) -- (node3SubPath1);
			\draw[->] (node3SubPath1) -- (node4SubPath1);
			\node[name=nodeSubPath1ComplementLabel,below=1cm of nodeSubPath1Label.east,anchor=east] {$q^c$};
			\node[name=node4SubPath1Complement,draw,circle] at (node4 |- nodeSubPath1ComplementLabel) {$4$};
			\node[name=node5SubPath1Complement,left=1cm of node4SubPath1Complement,draw,circle] {$5$};
			\node[name=node6SubPath1Complement,left=1cm of node5SubPath1Complement,draw,circle] {$6$};
			\draw[->] (node4SubPath1Complement) -- (node5SubPath1Complement);
			\draw[->] (node5SubPath1Complement) -- (node6SubPath1Complement);
			\node[name=nodeSubPath2Label,below=1cm of nodeSubPath1ComplementLabel.east,anchor=east] {$ p_{13} \in \subpath_{1 \bullet}(p_{14})$};
			\node[name=node1SubPath2,draw,circle] at (node1SubPath1 |- nodeSubPath2Label) {$1$};
			\node[name=node2SubPath2,left=1cm of node1SubPath2,draw,circle] {$2$};
			\node[name=node3SubPath2,left=1cm of node2SubPath2,draw,circle] {$3$};
			\draw[->] (node1SubPath2) -- (node2SubPath2);
			\draw[->] (node2SubPath2) -- (node3SubPath2);
			\node[name=nodeSubPath2ComplementLabel,below=1cm of nodeSubPath2Label.east,anchor=east] {$p_{34} $};
			\node[name=node3SubPath2Complement,draw,circle] at (node3SubPath2 |- nodeSubPath2ComplementLabel) {$3$};
			\node[name=node4SubPath2Complement,left=1cm of node3SubPath2Complement,draw,circle] {$4$};
			\draw[->] (node3SubPath2Complement) -- (node4SubPath2Complement);
		\end{tikzpicture}
	\end{center}
	\caption{An illustration of the idea of choosing the subpaths. {The
	complement of the subpath $ q^c $ is strictly shorter than the subpath $q$
	such that in the recursion \eqref{eqDefRecursionBStep} the left factor
	of the bracket has strictly smaller degree than the right factor{, hence \ref{itemPH4a} in the Definition of a P.~Hall basis holds}.}
	Also, the subpath $ p_{34} $ of the subpath $q$ is strictly shorter
	than $ q^c $ such that in the recursion \eqref{eqDefRecursionBStep} 
	the left factor of the right factor of the bracket has strictly
	smaller degree than the left factor of the bracket{, thus making sure that
	\ref{itemPH4b} holds as well}. 
	}
	\label{figIllustrationSubpaths}
\end{figure*}
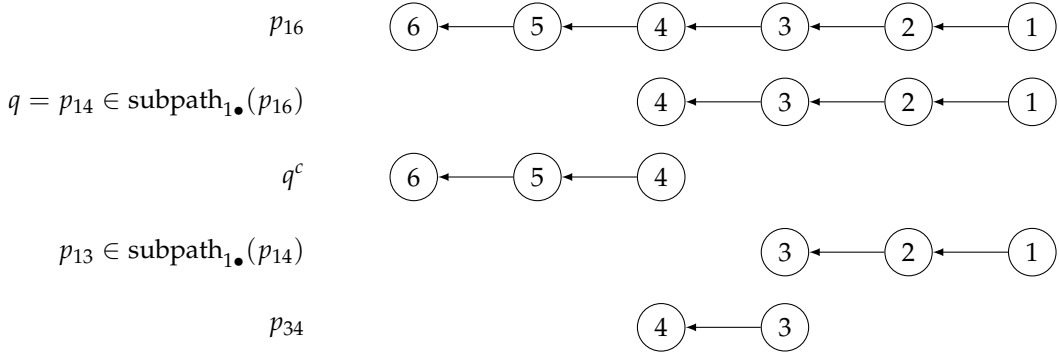
}{}
\begin{remark}\label{remarkChoiceStateToPassInformation}
\sout{It is worth pointing out, as become clear in the proof, that the aforementioned result is independent of the choice of $s$ as given in \eqref{eqChoiceS}; in fact, any
	{$ s \in \mathcal{I}\big(i_{\theta(p_{i_1 ir})}\big) = \lbrace i_{\theta(p_{i_1 i_r})}, n+ i_{\theta(p_{i_1 i_r})}, 2n + i_{\theta(p_{i_1 i_r})} \rbrace $}
	can be taken.
	Although the particular choice made does not make any difference in rewriting the
	non-admissible vector fields, it becomes relevant
	in designing suitable approximating inputs. 
	In particular, the choice of $s$ controls in which components of the complete
	state the perturbing inputs are injected.
	The specific choice \eqref{eqChoiceS}
	is motivated by the idea {of injecting} the most perturbation in the dual variables.
	Observe that the degrees of freedom for $s$ increase with the number 
	of constraints of each agent. In particular, it might as well
	happen that there is no degree of freedom if we do not augment
	the optimization problem \eqref{eqOptimizationProblem}.} 
	\oprocend
\end{remark}
{We no return to study~\eqref{eqSPArewritten}.}
Using \Cref{lemmaChoiceSubpath} we can then write~\eqref{eqSPArewritten}
as
\begin{align}
	\begin{split}
		\dot{\state}
		&=
		f_{\textup{adm}}(\state)
		-
		\sum\limits_{i=1}^n \sum\limits_{j=1}^n {\tilde{\MakeLowercase{\eqConstrL}}_{\textup{rest},ij}} \recBracketPHall{n+j,i}{}\big( p_{ij} \big) ( \state ) \\
		&-
		\sum\limits_{i=1}^n \sum\limits_{j=1}^n {\tilde{\MakeLowercase{\ineqConstrL}}_{\textup{rest},ij}} \recBracketPHall{2n+j,i}{}\big( p_{ij} \big) ( \state ) .
	\end{split}
	\label{eqSPArewrittenPHall}
\end{align}
and we can identify the set of brackets 
$ \mathcal{B} $ in \eqref{eqExtendedSystemLiu}
with
\begin{align}
	\begin{split}
	\mathcal{B} &= \big\lbrace \recBracketPHall{n+j,i}{}\big( p_{ij} \big): {\tilde{\MakeLowercase{\eqConstrL}}_{\textup{rest},ij}} \neq 0, i,j = 1,\dots,n \big\rbrace \\ 
	&\cup
	\big\lbrace \recBracketPHall{2n+j,i}{}\big( p_{ij} \big): {\tilde{\MakeLowercase{\ineqConstrL}}_{\textup{rest},ij}} \neq 0, i,j = 1,\dots,n  \big\rbrace,
	\end{split}
	\label{eqBracketsToBeExcited}
\end{align}
where now $ \mathcal{B} \subset \mathbb{P} $ for some P.~Hall
basis $ \mathcal{PH} = ( \mathbb{P}, \prec ) $,
and for the coefficients we have 
\begin{subequations}
\begin{align}
	v_{ \recBracketPHall{n+i,j}{}( p_{ij} )  }
	&= - {\tilde{\MakeLowercase{\eqConstrL}}_{\textup{rest},ij}} \textup{sign}\big(\recBracketPHall{n+i,j}{}( p_{ij} )(1)\big) \\
	v_{ \recBracketPHall{2n+i,j}{}( p_{ij} )  }
	&= - {\tilde{\MakeLowercase{\ineqConstrL}}_{\textup{rest},ij}} \textup{sign}\big(\recBracketPHall{n+i,j}{}( p_{ij} )(1)\big).
\end{align}
\end{subequations}
We are now ready to apply the algorithm presented in \cite{liu1997approximation}
to construct suitable approximating inputs.

\subsubsection{Approximating input sequences}\label{secApproximatingInputSequences}
We consider the collection of all agent dynamics \eqref{eqSPAAgentStructure}
given by
\begin{align}
	\dot{\state}^{\seqParam} = u^{\seqParam}(t,\state^{\seqParam}) 
	\ifthenelse{\boolean{longVersion}}{
	= 
	\begin{bmatrix}
		u_x^{\seqParam}(t,[x^{\seqParam},\dualEq^{\seqParam},\dualIneq^{\seqParam}]) \\
		u_{\dualEq}^{\seqParam}(t,[x^{\seqParam},\dualEq^{\seqParam},\dualIneq^{\seqParam}]) \\
		u_{\dualIneq}^{\seqParam}(t,[x^{\seqParam},\dualEq^{\seqParam},\dualIneq^{\seqParam}])
	\end{bmatrix}}{},
	\label{eqAllAgentsOriginalSystem}
\end{align}
where $ \state^{\seqParam} = [ {x^{\seqParam}}^\top, {\dualEq^{\seqParam}}^\top, {\dualIneq^{\seqParam}}^\top ]^\top $,
$ x^{\seqParam} \in \real^n $, {and} $ \dualEq^{\seqParam} \in \real^n $, $ \dualIneq^{\seqParam} \in \real^n $
are the stacked vectors of all $ x^{\seqParam}_i, \dualEq^{\seqParam}_i, \dualIneq^{\seqParam}_i $,
$ i = 1,2,\dots, n $, respectively, and 
\ifthenelse{\boolean{longVersion}}{
$ u_x^{\seqParam}, u_{\dualEq}^{\seqParam}, u_{\dualIneq}^{\seqParam}: \real \times \real^{3n} \to \real^n $
are the stacked vectors of all $ u_{x,i}^{\seqParam}, u_{\dualEq,i}^{\seqParam}, u_{\dualIneq,i}^{\seqParam} $,
$ i = 1,2,\dots,n $, respectively.}{
$u^\seqParam : \real \times \real^{3n} \to \real^{3n} $ is the stacked vector of
all $ u_{x,i}^{\seqParam}, u_{\dualEq,i}^{\seqParam}, u_{\dualIneq,i}^{\seqParam} $,
$ i = 1,2,\dots,n $.}
Following the algorithm presented in \cite{liu1997approximation}, we let the input take the form
\begin{align}
	u^{\seqParam}(t,\state^{\seqParam})
	&=
	f_{\textup{adm}}(\state^{\seqParam}) 
	+
	\sum\limits_{k = 1}^{M} \phi_k(\state^{\seqParam}) U_{k}^{\seqParam}(t),
	\label{eqDistributedControlInput}
\end{align}
where $ \Phi = \lbrace \phi_1, \phi_2,\dots,\phi_M \rbrace  $
is the set of admissible vector fields defined
in \eqref{eqDefSetOfAdmissibleVectorFields} and where $ \phi_k \in \mathbb{P} ${, $k~=~1,2,\dots,M$,}
for some P.~Hall basis $ \mathcal{PH}(\Phi) = ( \mathbb{P}, \prec ) $.
Further, $ U_{k}^{\seqParam}: \real \to \real $, $ k = 1,\dots, M$,
are so-called \emph{approximating input sequences} with sequence parameter
$ \seqParam \in {\natPos} $ which in the following we aim to construct 
in such a way that the solutions of \eqref{eqAllAgentsOriginalSystem}
uniformly converge to those of \eqref{eqSPArewrittenPHall} with increasing $ \seqParam $.
The algorithm in \cite{liu1997approximation} relies on {a
``superposition principle''}, i.e., we group all brackets in 
$ \mathcal{B} $ defined by \eqref{eqBracketsToBeExcited}
into equivalence classes, {which we later denote by $ E $}, {treat each equivalence class separately
and sum the resulting approximating inputs up in the end.
More precisely, we} associate to each class
an input $ U_{k,E}^{\sigma} $ and then let
\begin{align}
	U_{k}^{\seqParam}(t) = \sum\limits_{E \in \mathcal{E}} U_{k,E}^{\seqParam}(t), \label{eqInputSuperposition}
\end{align}
where $ \mathcal{E} $ is the set of all equivalence classes
in $ \mathcal{B} $. Roughly speaking, two brackets are said
to be equivalent if each vector field appears the same number
of times in the bracket but possibly in a different order.
A precise definition of the equivalence relation is given in
\Cref{defEquivalenceBrackets}.
For each {equivalence class} $ E \in \mathcal{E} $ and $ k = 1,\dots,M $ we then 
define the corresponding input $ U_{k,E}^{\seqParam}(t) $ as follows:
\begin{itemize}[leftmargin=*]
	\item If $ \degree_k(E) = 0 $: $ U_{k,E}^{\sigma}(t) = 0 $.
	\item If $ \degree(E) = 2, \degree_k(E) = 1 $:
	\begin{align}
		U_{k,E}^{\sigma}(t) = 2 \sqrt{\sigma} \Real \big( \eta_{E,k}(\omega_E) e^{\img \sigma \omega_{E} t} \big) . \label{eqInputBracket2}
	\end{align}
	\item If $ \degree(E) = N $, {$N \in \lbrace 3,4,\dots \rbrace $,} $ \degree_k(E) = 1 $:
	\begin{align}
		U_{k,E}^{\sigma}(t) = 2 \sigma^{\tfrac{N-1}{N}} \sum\limits_{\rho=1}^{\vert E \vert} \Real \big( \eta_{E}(\omega_{E,\rho,k}) e^{\img \sigma \omega t} \big). \label{eqInputBracketHigher}
	\end{align}
\end{itemize}
Here, it is $ \degree(E) = \degree(B) $, $ \degree_k(E) = \degree_k(B) $ for any $ B \in \mathcal{E} $.
Further, $ \omega_{E}, \omega_{E,\rho,k} \in \real $ are frequencies we will specify 
later, $ \eta_{E,k}, \eta_{E} : \real \to \mathbb{C} $ are coefficients to be chosen
{in dependence of} the frequencies, and $ \img \in \mathbb{C} $ is the imaginary unit.
{However, the superposition principle does not {hold} as desired {and} there are two major issues one has to take care of:}
\ifthenelse{\boolean{longVersion}}{
\begin{enumerate}[leftmargin=*]
    \item The input sequences $ U_{k,E}^{\seqParam} $ may not
   interfere with each other {in a way which ensures that}
    the superposition principle holds; 
    this can be dealt with by a proper choice of the 
    frequencies.
    \item Each input sequence $ U_{k,E}^{\seqParam} $ not only 
    {generates} 
    the desired brackets $E \cap \mathcal{B}$ {for $ \seqParam \to \infty $}, but also 
    all other equivalent brackets in $E$; 
    we can overcome this by a
    proper choice of the coefficients $ \eta_{\omega,k}, \eta_{\omega} $. 
   The idea behind {this is to also {generate}} 
    the undesired equivalent brackets {on purpose}, which itself
    also {generate} the desired brackets, in such a way that the
    undesired equivalent brackets all cancel out.
\end{enumerate}}{
    {(1) The input sequences $ U_{k,E}^{\seqParam} $ may not
   interfere with each other {in a way which ensures that}
    the superposition principle holds; 
    this can be dealt with by a proper choice of the 
    frequencies. 
    (2) Each input sequence $ U_{k,E}^{\seqParam} $ not only 
    {generates} 
    the desired brackets $E \cap \mathcal{B}$ {for $ \seqParam \to \infty $}, but also 
    all other equivalent brackets in $E$; 
    we can overcome this by a
    proper choice of the coefficients $ \eta_{\omega,k}, \eta_{\omega} $. 
   The idea behind {this is to also {generate}} 
    the undesired equivalent brackets {on purpose}, which itself
    also {generate} the desired brackets, in such a way that the
    undesired equivalent brackets all cancel out.} \\
}
While the problem at hand does not allow for simplifications in
the choice of the frequencies, the calculation of proper
coefficients $ \eta_{\omega,k}, \eta_{\omega} $ can be 
simplified drastically by exploiting some structural
properties of the set of brackets $ \mathcal{B} $. 
More precisely, there are two properties that turn out
to be beneficial: First, in each bracket $ B \in \mathcal{B} $
each vector field $ \phi_k $ appears only once, i.e., $ \delta_k(B) \in \lbrace 0,1 \rbrace $,
for any $ B \in \mathcal{B} $, $ k = 1,\dots,M $, and
second, for any bracket $ B \in \mathcal{B} $, all
equivalent brackets either evaluate to the same vector
field as $ B $ or vanish, see \ifthenelse{\boolean{longVersion}}{\Cref{lemmaEquivalentBrackets}}{\cite[Lemma 4]{mic2017distributedArxiv}}.
We present and discuss the simplified {calculation procedure} in \Cref{secAlgorithmApproxSeq}.
{While {this procedure}
may be tedious, it is not time-consuming, can be done off-line and is algorithmically
implementable.}
{It is worth mentioning that the calculation of the P.~Hall basis
as well as the approximating input sequences is not distributed and
requires preliminary global information; hence, the \emph{design} of the distributed algorithm
is not distributed but its implementation is. It is a 
matter of future research to develop distributed design procedures.}

\subsection{Distributed algorithm} 
We next state our main result which relates the solutions {of~\eqref{eqSPA}} with those of~\eqref{eqAllAgentsOriginalSystem}
in {closed loop} with the distributed control input \eqref{eqDistributedControlInput}-\eqref{eqInputBracketHigher}. {We use the notion of \emph{practically uniformly asymptotically stability} from~\cite{duerr2015thesis,duerr2013Lie}, without explicitly defining it here.}
\begin{theorem}\label{theoremMainResult}
Consider the distributed optimization problem \eqref{eqOptimizationProblem}
    and suppose that the communication topology is given by
    a strongly connected digraph with $n$ nodes. Assume that $ F $ is strictly convex and suppose
   further that \Cref{assConstraints1} - \ref{assConstraints2} hold. Consider the agent dynamics~\eqref{eqAllAgentsOriginalSystem} with the control law \eqref{eqDistributedControlInput}-\eqref{eqInputBracketHigher},
    {where the parameters in the control law are chosen according to
    the algorithm presented in~\Cref{secAlgorithmApproxSeq}.}
    Then, for each $ \varepsilon > 0 $, {for each $ T > 0 $,} and for each initial condition
    $ \state^{\sigma}(0) = \state(0) = \state_0 \in \mathcal{R}(\mathcal{M}) $, {with $ \mathcal{R}(\mathcal{M}) $ given in~\eqref{eqDefRegionOfAttraction},} there exists $ \seqParam^* > 0 $
    such that for all $ \seqParam > \seqParam^* $ the following holds:
   For all {$ 0 \leq t \leq T $}, {we have}
    \begin{align}
        \big\Vert \state^{\seqParam}(t) - \state(t) \big\Vert \leq \varepsilon, \label{eqLemmaStabResultTrajApprox}
    \end{align} 
    where $ \state^{\seqParam}(t) $ is the solution {of}~\eqref{eqAllAgentsOriginalSystem}
    with the control law \eqref{eqDistributedControlInput} - \eqref{eqInputBracketHigher} and
    $ \state(t) = \big( x(t),\dualEq(t),\dualIneq(t) \big) $ is the solution of
    \eqref{eqSPA}, with initial condition $ \state^{\sigma}(0) = \state(0) = \state_0 $. 
    Further, the set $ \mathcal{M} $ defined by~\eqref{eqDefSetOfSaddlePoints}
    is practically uniformly asymptotically stable {(given $\degree(B)=2$ for all $ B\in\mathcal{B}$ in~\eqref{eqBracketsToBeExcited})}. \oprocend
\end{theorem}
{We postpone the proof of this result to \Cref{secAppendixProofMainResult} and focus on its useful implications in the next section.}

\ifthenelse{\boolean{longVersion}}{
\subsection{Filtered saddle-point dynamics}\label{secFilteredSPD}
The highly oscillatory nature of the approximating inputs naturally
leads to an undesired oscillating behavior of the closed-loop
trajectories of the distributed approximation. 
As discussed in~\Cref{secAlgorithmApproxSeq}, the effect
on the primal variables, which are in most cases the ones one 
is most interested in, can be reduced by a proper design of the
approximating inputs. Another natural remedy to this problem
is to make use of filters which we want to briefly discuss
in the following. There are different ways of introducing 
filters in the feedback loop; in the following we concentrate
on the situation depicted in~\Cref{figFilteredSPD}, where only
the signal $ u_x, u_\dualEq, u_\dualIneq $ are modified
by means of low-pass filters $ G_x, G_\dualEq, G_\dualIneq $,
where $ G_x, G_\dualEq, G_\dualIneq $ are square stable and proper
transfer matrices of appropriate dimension. In view of a distributed
implementation we restrict ourselves to diagonal transfer matrices;
hence the additional filters do not introduce new variables
which are not available to an agent in a distributed setting.
These filtered saddle-point dynamics can also be interpreted as
\emph{higher order saddle-point dynamics} where the minimization in the 
primal variable as well as the maximization in the dual variables
is not performed by means of a standard gradient descent or ascent, 
respectively, but higher order optimization algorithms~\cite{mic2014heavy}
are used. A thorough analysis of these filtered saddle-point dynamics
is still open, but we emphasize that, as long as the filters are
``sufficiently fast'', similar stability results can be obtained
making use of singular perturbation theory. 

As to the distributed approximation of the filtered saddle-point
dynamics, only minor modifications are required. In rough words,
the non-admissible terms appearing in the filtered
saddle-point dynamics take the same form as the ones without a filter 
but, since the complete state is augmented by the internal states
of the filter, they appear in a different component.
Hence, we basically only need to adapt the index sets~\eqref{eqDefIndexSetAgent}
and augment the vector fields~\eqref{eqDefAdmissibleVectorFields}.
We illustrate the effect of additional filters by means
of an example in~\Cref{secSimExample}.

\begin{figure}[ht]
	\begin{center}
		\begin{tikzpicture}[>=latex]
			\node[name=nodeAgentDynamics,draw,rectangle,inner sep=5pt,minimum width=1.9cm,minimum height=1.9cm] {
			\begin{minipage}{1.3cm}\vspace*{-1.1em}
				\begin{align*}
					\dot{x} &= u_{x,\textup{fil}} \\
					\dot{\dualEq} &= u_{\dualEq,\textup{fil}} \\
					\dot{\dualIneq} &= u_{\dualIneq,\textup{fil}}
				\end{align*}
			\end{minipage}
			};
			\node[name=nodeController,below=0.7cm of nodeAgentDynamics.south east,anchor=north east,draw,rectangle,inner sep=5pt,minimum width=2cm,minimum height=2cm] {
			\begin{minipage}{1.3cm}\vspace*{-1.1em}
				\begin{align*}
					u_\dualIneq &= \textup{diag}(\dualIneq) \nabla_{\dualIneq} L(x,\dualEq,\dualIneq) \\
					u_\dualEq &= \nabla_{\dualEq} L(x,\dualEq,\dualIneq) + w(\dualEq) \\
					u_x &= \nabla_x L(x,\dualEq,\dualIneq) 
				\end{align*}
			\end{minipage}
			};
			\node[name=nodeFilterPrimal,draw,rectangle,anchor=east,minimum width=1.2cm,yshift=19pt,anchor=west] at (nodeController.west |- nodeAgentDynamics.west) {$ G_x(s) $};
			\node[name=nodeFilterDualEq,draw,rectangle,anchor=east,minimum width=1.2cm,yshift=0pt,anchor=west] at (nodeController.west |- nodeAgentDynamics.west) {$ G_\dualEq(s) $};
			\node[name=nodeFilterDualIneq,draw,rectangle,anchor=east,minimum width=1.2cm,yshift=-19pt,anchor=west] at (nodeController.west |- nodeAgentDynamics.west) {$ G_\dualIneq(s) $};
			\draw[->] (nodeFilterPrimal) -- node[anchor=south,inner sep=2pt] {$u_{x,\textup{fil}}$} (nodeAgentDynamics.west |- nodeFilterPrimal.east);
			\draw[->] (nodeFilterDualEq) -- node[anchor=south,inner sep=2pt] {$u_{\dualEq,\textup{fil}}$} (nodeAgentDynamics.west |- nodeFilterDualEq.east);
			\draw[->] (nodeFilterDualIneq) -- node[anchor=south,inner sep=2pt] {$u_{\dualIneq,\textup{fil}}$} (nodeAgentDynamics.west |- nodeFilterDualIneq.east);
			\draw[->] (nodeAgentDynamics.east) -- ++(1cm,0) |- (nodeController.east);
			\draw[->] ([yshift=16pt] nodeController.west) -- ++(-0.5cm,0) |- node[pos=1,anchor=south east,inner sep=2pt] {$u_{\dualIneq\vphantom{,}}$} (nodeFilterDualIneq);
			\draw[->] ([yshift=0pt] nodeController.west) -- ++(-0.75cm,0) |- node[pos=1,anchor=south east,inner sep=2pt] {$u_{\dualEq\vphantom{,}}$} (nodeFilterDualEq);
			\draw[->] ([yshift=-16pt] nodeController.west) -- ++(-1cm,0) |- node[pos=1,anchor=south east,inner sep=2pt] {$u_{x\vphantom{,}}$} (nodeFilterPrimal);
		\end{tikzpicture}
	\end{center}
	\caption{Saddle-point dynamics~\eqref{eqSPA} with additional low-pass filters $ G_x, G_\dualEq, G_\dualIneq $.}\label{figFilteredSPD}
\end{figure}

}{}

\section{Special cases and examples}
In this section we discuss special cases in which the inputs can be
given in explicit form and present several simulation examples
illustrating the previous results.

\subsection{Explicit representation of approximating inputs for low order brackets}\label{secExplicitFormulaDegreeThree}
While the algorithm given in \Cref{secAlgorithmApproxSeq} can in general
be complicated to implement,
{{the} procedure becomes particularly simple to implement 
in scenarios} where the set of brackets
$ \mathcal{B} $ defined in~\eqref{eqBracketsToBeExcited} only contains
brackets of degree less or equal than three.
As stated {in our next result,} in this case the set of equivalent brackets
only contains the bracket itself but no other bracket, thus the second
issue 2 in \Cref{secApproximatingInputSequences} {does not come into play}. 
\begin{proposition}\label{lemmaLowDegreeBrackets}
Consider \eqref{eqSPArewrittenPHall} and assume that all paths $ p_{ij} $
	fulfill $ \length(p_{ij}) \leq 3 $. Let $ \mathcal{PH}(\Phi) = ( \mathbb{P}, \prec ) $
	be any P.~Hall basis of $ \Phi $ defined by \eqref{eqDefSetOfAdmissibleVectorFields}
	{that fulfills $ h_{k_1,k_2} \prec h_{k_3,k_4} $} for all $ k_4 > k_2 $. 
	Then, for any path $ p_{ij} $ with $ \length(p_{ij}) \leq 3 $,
	we have that {the equivalence class corresponding to
	the bracket $ \recBracketPHall{r+j,i}{}( p_{ij} ) $ fulfills}
	\begin{align}
		E_{ \recBracketPHall{r+j,i}{}( p_{ij} ) }
		&= \lbrace B \in \mathbb{P} : B \sim \recBracketPHall{r+j,i}{}( p_{ij} ), B(z) \not\equiv 0 \rbrace \nonumber \\
		&= \lbrace \recBracketPHall{r+j,i}{}( p_{ij} ) \rbrace
		\label{eqLemmaLowDegree}
	\end{align}
	for $ r \in \lbrace n,2n \rbrace $, where the equivalence relation $ \sim $ is defined
	by \Cref{defEquivalenceBrackets}.
	\oprocend
\end{proposition}
\begin{remark}
It should be noted that the ordering of the P.~Hall basis is important
	for this result to hold. 
	Further, {if \Cref{assConstraints2} does not hold,}
	different brackets are introduced in \eqref{eqSPArewrittenPHall} 
	which still are of degree three under the assumption that all paths $ p_{ij} $
	fulfill $ \length(p_{ij}) \leq 3 $ but have a different structure.
	Hence, the assumption on the ordering is {in general} not sufficient anymore. 
	\oprocend
\end{remark}
A proof of this result can be found in \Cref{secProofLemmaLowDegreeBrackets}.
The condition that all paths $p_{ij}$ in \eqref{eqSPArewrittenPHall} are of 
length less or equal than three holds, for example, if the longest cordless cycle in
$ \mathcal{G} $ is of length $4$.
Using the result of \Cref{lemmaLowDegreeBrackets} and following the algorithm presented in \Cref{secAlgorithmApproxSeq}, we
obtain \\[0.5em]
\labelitemi\hspace{2pt} if $ {E = \lbrace B \rbrace = \lbrace [ \phi_{k_1}, \phi_{k_2} ] \rbrace } $: \\[-2em]
	\begin{align}
		&U_{k,E}^{\seqParam}(t) \\
		&=
		\begin{cases}
			\hfill - \sqrt{2 \seqParam} \tfrac{1}{\beta_E} \sqrt{ \vert v_B \omega_E \vert } \cos(\seqParam\omega_{{E}} t)         & \text{if } k = k_1 \\
			\sign(v_B\omega_B) \sqrt{2 \seqParam} \beta_E \sqrt{\vert v_B \omega_E \vert} \sin(\seqParam\omega_{{E}} t) & \text{if } k = k_2 \\
			\hfill 0                                                                                               & \text{otherwise,}
		\end{cases} \nonumber
	\end{align}
\labelitemi\hspace{2pt} if {$ E = \lbrace B \rbrace = \lbrace \big[ \phi_{k_1}, [\phi_{k_2}, \phi_{k_3}] \big] \rbrace $}: 
	\begin{align}
		&U_{k,E}^{\seqParam}(t) \\
		&=
		\begin{cases}
			\hfill - \seqParam^{\tfrac{2}{3}} 2 \beta_E ( \omega_{E,k_1} \omega_{E,k_2} )^{\tfrac{1}{3}} \cos(\seqParam \omega_{E,k} t)        & \text{if } k = k_1, k_3 \\
			- \seqParam^{\tfrac{2}{3}} 2 \tfrac{1}{\beta_E^2} ( \omega_{E,k_1} \omega_{E,k_2} )^{\tfrac{1}{3}} \cos(\seqParam \omega_{E,k_2} t) & \text{if } k = k_2 \\
			\hfill 0                                                                                                                     & \text{otherwise,}
		\end{cases} \nonumber
	\end{align}
where $ \beta_E \neq 0 $ is a design parameter. The frequencies
$ \omega_E, \omega_{E,k} \in \real \setminus \lbrace 0 \rbrace $
need to be chosen such that they fulfill the following properties: 
\begin{itemize}[leftmargin=*]
	\item All frequencies $ \omega_E $, $ E \in \mathcal{E} $, $ \degree(E) = 2 $, are distinct.
	\item For each $ E = {\lbrace} B {\rbrace} = \big[ \phi_{k_1}, [\phi_{k_2}, \phi_{k_3}] \big] $,
	the set of frequencies $ \lbrace \omega_{E,k_1}, \omega_{E,k_2}, \omega_{E,k_3} \rbrace $
	is minimally canceling, see \Cref{defSymmetricallyMinimallyCanceling}.
	\item The collection of sets \\[-2em]
	\begin{align*}
		\big\lbrace \lbrace \omega_E \rbrace_{E \in \mathcal{E}, \degree(E) = 2}, \lbrace \omega_{E,k_1}, \omega_{E,k_2}, \omega_{E,k_3} \rbrace_{E \in \mathcal{E}, \degree(E) = 3} \big\rbrace
	\end{align*}\\[-2em]
	is an independent collection, see \Cref{defIndependentCollection}.
\end{itemize}
{Note that there always exist frequencies that fulfill these
properties, see~\cite{liu1997approximation}.}
Similar explicit formulas can as well
be obtained for brackets of higher degree but they become more complicated. The main reason is that, while for brackets
of degree strictly less than four all equivalent brackets evaluate to
zero\ifthenelse{\boolean{longVersion}}{ (cf. \Cref{tableNonZeroBrackets})}{}, this is no longer the case for brackets of higher degree such
that now the second issue discussed in \Cref{secApproximatingInputSequences}
needs to be taken care of.

\subsection{Simulation examples}\label{secSimExample}
\begingroup
\arraycolsep=3pt\def\arraystretch{1}
\begin{figure*}[ht]
    \def\deltaAngle{20}
	\centering
	\begin{minipage}[t]{0.47\textwidth}\vspace{0pt}
		\begin{tikzpicture}[>=latex]
			\node at (90:1.7) {Graph (a)};
			\node[circle,draw,name=node1,anchor=center] at(90+36:1.3) {$1$};
			\node[circle,draw,name=node2,anchor=center] at(90-36:1.3) {$2$};
			\node[circle,draw,name=node3,anchor=center] at(90-36-72:1.3) {$3$};
			\node[circle,draw,name=node4,anchor=center] at(90-36-72-72:1.3) {$4$};
			\node[circle,draw,name=node5,anchor=center] at(90-36-72-72-72:1.3) {$5$};
			\draw[->] (node1) -- (node2);
			\draw[->] (node1) -- (node5);
			\draw[->] (node2) -- (node3);
			\draw[->] (node3) -- (node1);
			\draw[->] (node4) -- (node3);
			\draw[->] (node5) -- (node4);
			\draw[->] (node5) -- (node2);
			\draw[->,fictitiousColor,dashed] (node2) to[out=180-\deltaAngle,in=\deltaAngle] (node1); 
			\draw[->,fictitiousColor,dashed] (node3) to[out=108-\deltaAngle,in=-72+\deltaAngle] (node2);
			\draw[<-,fictitiousColor,dashed] (node4) to[out=36-\deltaAngle,in=-144+\deltaAngle] (node3);
			\draw[->,fictitiousColor,dashed] (node2) to[out=-144+\deltaAngle,in=36-\deltaAngle] (node5);
		\end{tikzpicture}
		\hspace*{1em}
		\begin{tikzpicture}[>=latex]	
			\node at (90:1.7) {Graph (b)};
			\node[circle,draw,name=node1,anchor=center] at(90+36:1.3) {$1$};
			\node[circle,draw,name=node2,anchor=center] at(90-36:1.3) {$2$};
			\node[circle,draw,name=node3,anchor=center] at(90-36-72:1.3) {$3$};
			\node[circle,draw,name=node4,anchor=center] at(90-36-72-72:1.3) {$4$};
			\node[circle,draw,name=node5,anchor=center] at(90-36-72-72-72:1.3) {$5$};
			\draw[->] (node1) -- (node2);
			\draw[->] (node1) -- (node5);
			\draw[->] (node2) -- (node3);
			\draw[->] (node3) -- (node1);
			\draw[->] (node4) -- (node3);
			\draw[->] (node5) -- (node4);
			\draw[->,fictitiousColor,dashed] (node2) to[out=180-\deltaAngle,in=\deltaAngle] (node1); 
			\draw[->,fictitiousColor,dashed] (node3) to[out=108-\deltaAngle,in=-72+\deltaAngle] (node2);
			\draw[<-,fictitiousColor,dashed] (node4) to[out=36-\deltaAngle,in=-144+\deltaAngle] (node3);
			\draw[->,fictitiousColor,dashed] (node2) to[out=-144+\deltaAngle,in=36-\deltaAngle] (node5);
			\draw[<-,fictitiousColor,dashed] (node2) to[out=-144-\deltaAngle,in=36+\deltaAngle] (node5);
		\end{tikzpicture}
    \end{minipage}
    \begin{minipage}[t]{0.51\textwidth}\vspace*{1em}
	\ifthenelse{\boolean{longVersion}}{\setlength{\tabcolsep}{2pt}}{\setlength{\tabcolsep}{4pt}}
	\begin{tabular}{@{}lll@{}}
		\toprule 
		vector         & corresponding                                                       & Lie bracket  \\ 
		field          & path                                                                & representation \\ \midrule 
		$ h_{n+2,3} $  & $ \pathLeft v_3 \pathSep v_1 \pathSep v_2 \pathRight $              & $ \big[ h_{n+2,n+1}, h_{n+1,3} \big] $ \\[0.2em]
		$ h_{2n+1,2} $ & $ \pathLeft v_1 \pathSep v_3 \pathSep v_1 \pathRight $              & $ \big[ h_{2n+1,2n+3}, h_{2n+3,3} \big] $ \\[0.2em]
		$ h_{n+5,2} $  & $ \pathLeft v_2 \pathSep v_3 \pathSep v_1 \pathSep v_5 \pathRight $ & $\big[ h_{n+5,n+1}, [ h_{n+1,n+3}, h_{n+3,2} ] \big] $ \\[0.2em]
		$ h_{2n+4,3} $ & $ \pathLeft v_3 \pathSep v_1 \pathSep v_5 \pathSep v_4 \pathRight $ & $\big[ h_{2n+4,2n+5}, [ h_{2n+5,2n+1}, h_{2n+1,3} ] \big] $ \\
		\bottomrule
	\end{tabular}
	\end{minipage}
	\caption{Left: Two communication graphs~(a) and (b) for the simulation example from \Cref{secSimExample}.
	The dashed green arrows indicate the required fictitious edges, respectively.
	Right: The results of applying \Cref{lemmaConstructionBrackets} to rewrite the non-admissible
	vector fields in the example from \Cref{secSimExample} in terms of Lie brackets
	of admissible vector fields ($n=5$).}
	\label{figGraphExample}
\end{figure*}
Next, we present some {simulated} examples to illustrate our results: 
{We} consider an optimization problem of the form \eqref{eqOptimizationProblem}
with $ n = 5 $ agents, where, for $ i = 1,2,\dots, 5 $,
$
	F_i(x_i) = (x_i-i)^2,
$
and the constraints are given by
\begin{subequations}
\begin{align}
	x_1 - x_2 {\leq} -10, \quad
	x_2 - x_3 = 1, \\
	x_4 + x_3 {\leq} -3, \quad
	x_5 - x_2 = 7, 
\end{align}\label{eqExampleConstraints}
\end{subequations}
such that after augmentation we have for the matrices that define
the constraints in \eqref{eqOptimizationProblemAugmented}
\begin{align}
	\eqConstrL &= 
	\left[
	\begin{array}{rrrrr}
	0 &  0 &  0 & 0 & 0 \\
	0 &  1 & -1 & 0 & 0 \\
	0 &  0 &  0 & 0 & 0 \\
	0 & \hphantom{-} 0 & \hphantom{-}0 & \hphantom{-}0 & \hphantom{-}0 \\
	0 & -1 &  0 & 0 & 1
	\end{array}\right], &
	\eqConstrR &= \left[\begin{array}{l} 0 \\ 1 \\ 0 \\ 0 \\ 7 \end{array}\right], \\
	\ineqConstrL &= 
	\left[
	\begin{array}{rrrrr}
	1 & -1 & 0 & 0 & 0 \\
	0 &  0 & 0 & 0 & 0 \\
	0 &  0 & 0 & 0 & 0 \\
	0 &  0 & 1 & 1 & 0 \\
	0 & \hphantom{-}  0 & \hphantom{-}0 & \hphantom{-}0 & \hphantom{-}0
	\end{array}\right], &
	\ineqConstrR &= \left[\begin{array}{r} -10 \\ K \\ K \\ -3 \\ K \end{array}\right], 
\end{align}
where $ K = 3 $ but can as well be chosen arbitrary as long as $ K > 0 $.
We consider two different communication graphs as depicted in
\Cref{figGraphExample}, where graph~(b) is the same as {graph~(a)} 
except that the {edge} from agent 5 to agent 2 got broken,
thus an additional fictitious edge is required.
While the constraints match the communication topology of graph~(a),
i.e., \Cref{assConstraints2} holds, this is not the case for graph~(b)
due to the last constraint in \eqref{eqExampleConstraints}.
We first consider the case that graph~(a) represents the communication
topology. In this case, the graph Laplacian is given by
\begin{align}
	\laplacian
	=
	\left[
	\begin{array}{rrrrr}
	2  & -1 &  0 &  0 & -1 \\
	0  &  1 & -1 &  0 &  0 \\
	-1 &  0 &  1 &  0 &  0 \\
	0  &  0 & -1 &  1 &  0 \\
	0  & -1 &  0 & -1 &  2
	\end{array}
	\right]
\end{align}
and hence
\begin{align}
	\tilde{\eqConstrL}_{\textup{adm}} =
	\begin{bmatrix}
	0 & 0 & 0 & 0 & 0 \\
	0 & 1 & 0 & 0 & 0 \\
	0 & 0 & 0 & 0 & 0 \\
	0 & 0 & 0 & 0 & 0 \\
	0 & 0 & 0 & 0 & 1
	\end{bmatrix}, 
	\tilde{\ineqConstrL}_{\textup{adm}} =
	\begin{bmatrix}
	1 & 0 & 0 & 0 & 0 \\
	0 & 0 & 0 & 0 & 0 \\
	0 & 0 & 0 & 0 & 0 \\
	0 & 0 & 0 & 1 & 0 \\
	0 & 0 & 0 & 0 & 0
	\end{bmatrix}.
\end{align}
The saddle-point dynamics \eqref{eqSPASeparated} are then given by
\begin{align}
	\dot{\state} 
	&= f_{\textup{adm}}(\state) - ( - \unitVec{3} \state_7 -  \unitVec{2} \state_{11} - \unitVec{2} \state_{10} + \unitVec{3} \state_{14} ) \nonumber \\
	&= f_{\textup{adm}}(\state) + h_{n+2,3}(\state) + h_{2n+1,2}(\state)  \\
	&\phantom{= f_{\textup{adm}}(\state)}\; + h_{n+5,2}(\state) - h_{2n+4,3}(\state), \nonumber
\end{align} 
where the admissible part $ f_{\textup{adm}}: \real^{15} \to \real^{15} $ is defined by \eqref{eqDefAdmissiblePart}
and the remaining four vector fields are non-admissible.
Following \Cref{lemmaConstructionBrackets} and choosing the subpaths as suggested
in \Cref{lemmaChoiceSubpath} we then rewrite the non-admissible vector fields as given in
{the table in~\Cref{figGraphExample}.}

As a next step, we need to {write the Lie brackets as a linear combination of brackets in some}
P.~Hall basis $ \mathcal{PH}(\Phi) = ( \mathbb{P}, \prec ) $ with
\begin{align*}
	\Phi \, = \, \lbrace \, h_{n+2,n+1}, \, h_{n+1,3}, \, h_{2n+1,2n+3}, \, h_{2n+3,3}, \, h_{n+5,n+1}, \\
	h_{n+1,n+3}, \, h_{n+3,2}, \, h_{2n+4,2n+5}, \, h_{2n+5,2n+1}, \, h_{2n+1,3} \rbrace.
\end{align*}
In general, we can choose any P.~Hall basis and then make 
use of \Cref{remarkProjectionPHall} for the projection. 
However, in this case it is also easily possible to properly choose
the ordering of the P.~Hall basis in such a way that the brackets
in {the table in~\Cref{figGraphExample}} are already in $ \mathbb{P} $.
More precisely, we only have to make sure that
$ h_{n+2,n+1} \prec h_{n+1,3} $, $ h_{2n+1,2n+3} \prec h_{2n+3,3} $,
$ h_{n+1,n+3} \prec h_{n+3,2} $, $ h_{n+1,n+3} \prec h_{n+5,n+1} $,
$ h_{2n+5,2n+1} \prec h_{2n+1,3} $, $ h_{2n+5,2n+1} \prec h_{2n+4,2n+5} $.
Note that this is in general not possible, since the conditions
might be conflicting and -- to keep this example more general --
we do not {adapt the ordering in that way} in our implementation.

We are now ready to apply the algorithm presented in \Cref{secAlgorithmApproxSeq}.
We do not discuss the resulting input sequences in detail here
and also do not provide the {complete} simulation results due to space limitations,
but instead do this for the case that the communication graph is given by
graph~(b). We refer the interested reader to employ the Matlab implementation
{provided in the supplementary material\ifthenelse{\boolean{longVersion}}{}{~in~\cite{mic2017distributedArxiv}}.}
We next discuss the implications of having the communication
graph given by graph~(b) in \Cref{figGraphExample} instead
of graph~(a). {Since the {edge} from node~2 to node~3 is missing in the graph,} \Cref{assConstraints2}
does no longer hold.
{In particular,} the vector field $ h_{n+2,n+5}(\state) = \state_{n+2} \unitVec{n+5} $,
which is included in the admissible vector field $ f_{\textup{adm}} $ in
case the communication is given by graph~(a), now is non-admissible. Despite \Cref{assConstraints2} not being fulfilled,
we can {still} use \Cref{lemmaConstructionBrackets} to rewrite
$ h_{2,n+5} $, since the result is completely independent of this assumption. {Indeed}, the corresponding path is given
by $ p_{52} = \pathLeft v_5 \pathSep v_4 \pathSep v_3 \pathSep v_1 \pathSep v_2 \pathRight $
and we obtain
\begin{align}
	&h_{2,n+5}(\state) \\
	&= \big[ [ h_{2,n+1}, h_{n+1,n+3} ], [ h_{n+3,n+4}, h_{n+4,n+5} ] \big](\state). \nonumber
\end{align}
We can then follow the same procedure as discussed before to project
on any P.~Hall basis, where $ \Phi $ now additionally includes
the vector fields $ h_{2,n+1} $, $h_{n+1,n+3}$, $h_{n+3,n+4}$, and $h_{n+4,n+5} $,
and then apply the algorithm presented in \Cref{secAlgorithmApproxSeq}.
The {corresponding} simulation results are depicted in \Cref{figSimExampleGraphB}.
\ifthenelse{\boolean{longVersion}}{
{
As already indicated in~\Cref{secFilteredSPD}, the nature of the approximating
inputs produces heavy oscillations in the agents' states. We next want to illustrate
how a properly chosen filter as described in~\Cref{secFilteredSPD} can be
used to dampen these oscillations while still maintaining the distributed structure.
In this example, we assume that the low-pass filters $ G_x, G_\dualEq, G_\dualIneq $
are of first order and take $ G_x(s) = \tfrac{70}{s+70} $, $ G_\dualEq(s) = \tfrac{170}{s+170} $,
$ G_\dualIneq(s) = \tfrac{170}{s+170} $. We do not go through the calculations
necessary to find distributed approximations of the filtered saddle-points dynamics
since they literally follow the lines of the first part of the example.
The corresponding simulation results are depicted in~\Cref{figSimExampleGraphB}. Compared to~\Cref{figSimExampleGraphB},
the trajectories of the distributed approximation show less oscillations
which in turn also leads to an improved approximate solution of the 
distributed optimization problem.}
}{{As apparent from the numerical results, the nature of the approximating
inputs produces oscillations in the agents' states. By adding properly
designed filters, their amplitude can be reduced; we discuss this in more 
detail in~\cite{mic2017distributedArxiv}}}.

\begin{figure*}[htp]
	\begin{center}
		\ifthenelse{\boolean{longVersion}}{
		\includegraphics[width=0.8\textwidth]{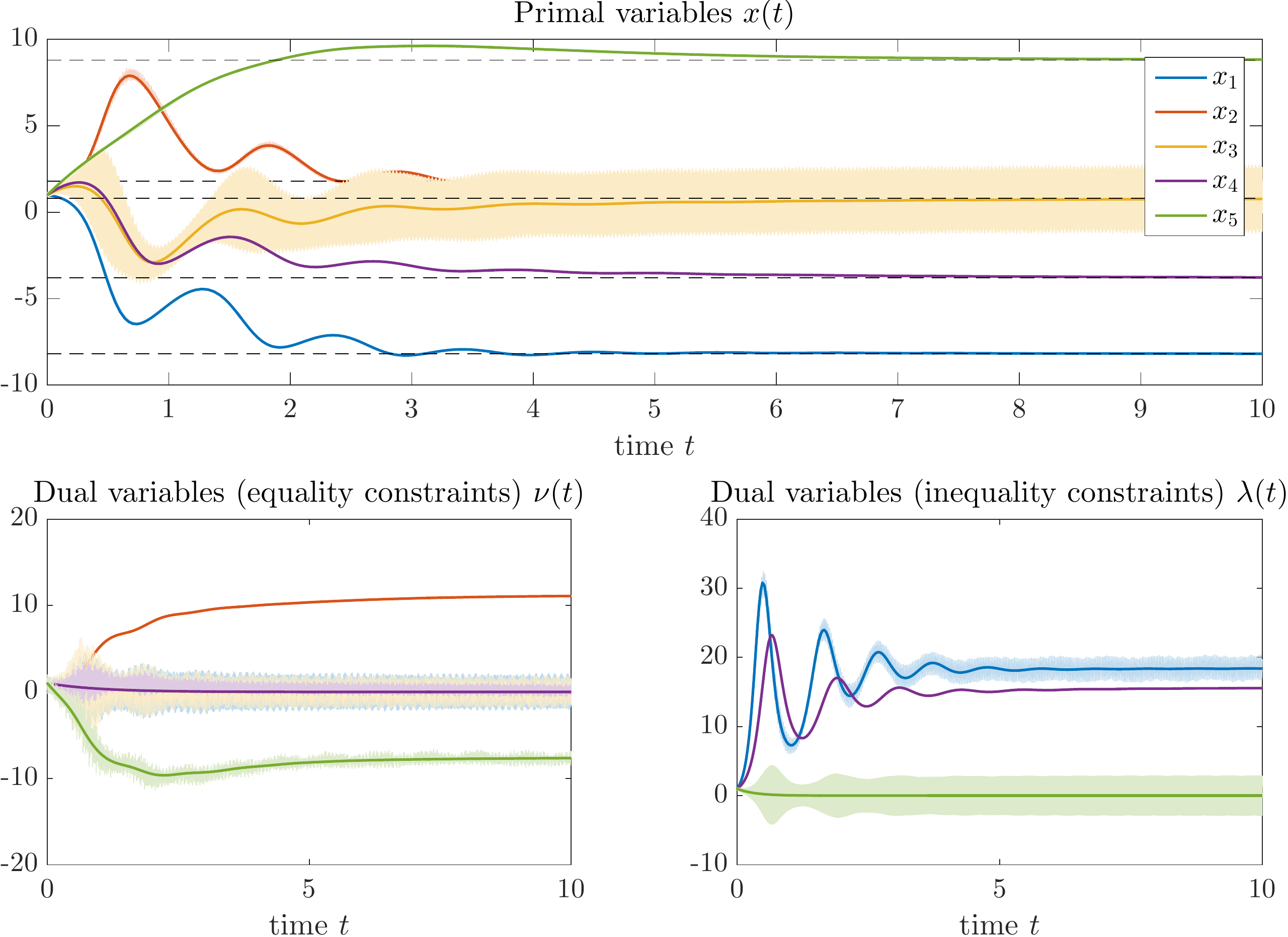} 
		\includegraphics[width=0.8\textwidth]{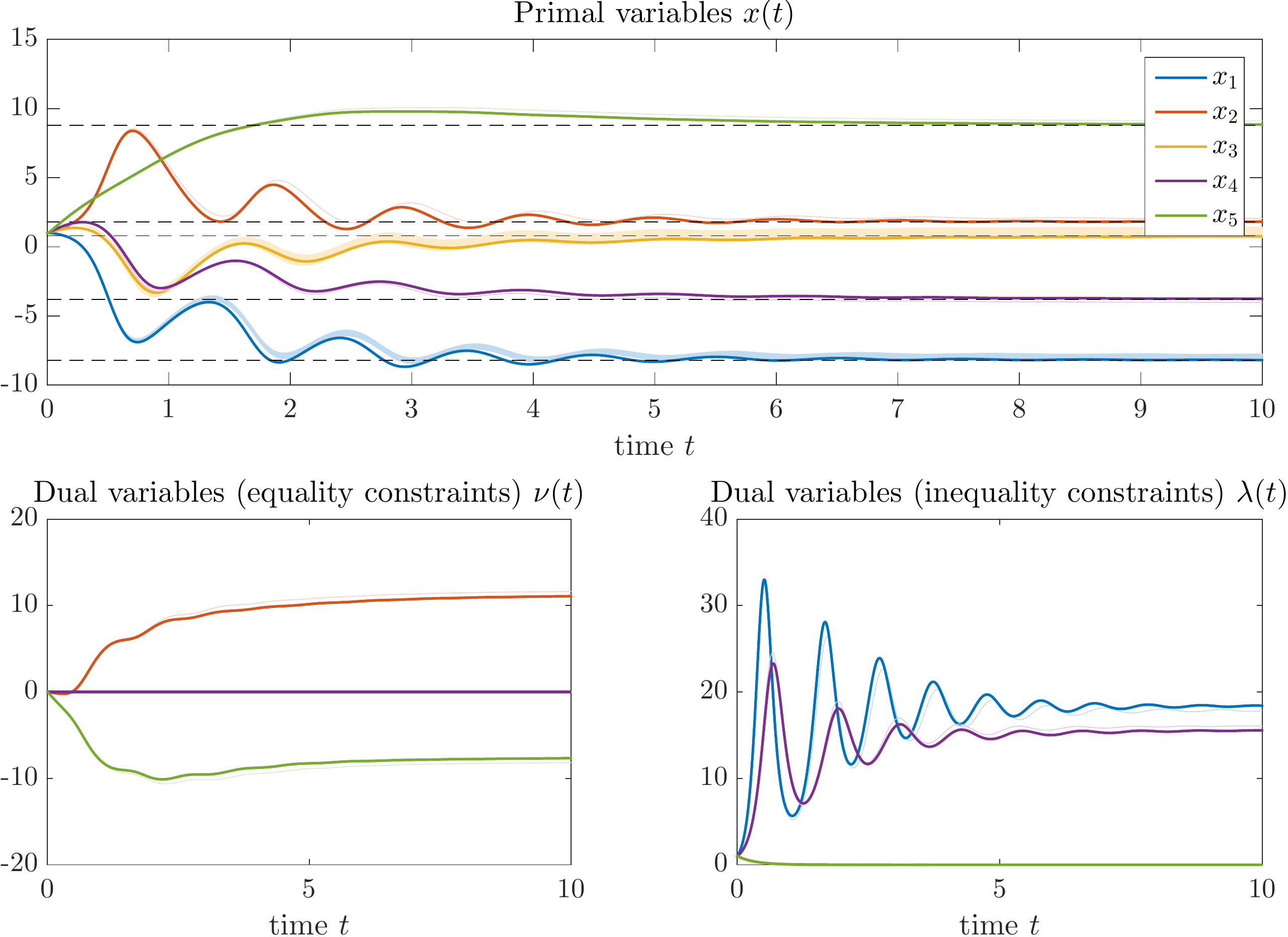} 
		}{\includegraphics[width=0.96\textwidth]{figs/simExampleGraphB.pdf} }
	\end{center}
	\caption[]{Simulation results for the example {of} \Cref{secSimExample} with
	communication graph~(b) {given in} \Cref{figGraphExample}\ifthenelse{\boolean{longVersion}}{ without (top) and with additional filters (bottom)}{}.
	The thick lines depict
	the trajectories of the (non-distributed) saddle-point dynamics with initial
	condition $ \state(0) = \mathbf{1} \in \real^{15} $, whereas the thinner oscillating
	lines depict the solution of the distributed approximation with the same
	initial condition $ \state^{\seqParam}(0) = \state(0) $. Where no oscillating
	lines are visible, they are covered by the corresponding component of the solution $ z(\cdot) $.
	The dashed black lines indicate the optimal solution of the optimization problem
	given by $ x^\star = [ -8.2, \; 1.8, \; 0.8, \; -3.8, \; 8.8 ]^\top $.
	\ifthenelse{\boolean{longVersion}}{For both simulations the frequencies
	were chosen differently but according to some heuristics making sure that the 
	minimally canceling property {from~\Cref{defSymmetricallyMinimallyCanceling}}
	is fulfilled.}{The frequencies required in the approximating inputs were chosen according to some heuristics
	with absolute values in the range between $ 2.6309 $ and $ 62.6381 $
	and making sure that the minimally canceling property {from~\Cref{defSymmetricallyMinimallyCanceling}} is fulfilled.}
	Further, {we used} $ \sigma = 1000$.}
	\label{figSimExampleGraphB}
\end{figure*}
\endgroup

\section{Conclusion and outlook}
We presented a new approach to distributed optimization problems
where the communication topology is given by a directed graph. 
Our approach is based on a two-step procedure where in a first
step first we derived suitable Lie bracket representations of
saddle-point dynamics and then used Lie bracket approximations 
techniques from geometric control theory to obtain 
distributed control laws. While we limited ourselves to 
the class of convex problems with separable cost function
and linear equality and inequality constraints that match the communication
topology, the methodology is applicable
to a much larger class of optimization problems including,
for example, non-linear constraints, {constraints not compatible with the graph structure}
or non-separable cost functions{; we discuss in~\cite{mic2018extensions}
how the rewriting procedure has to be adapted}. 
{Certainly, this generality comes
with the cost of a possibly complex calculation of the approximating inputs;
however, the strength of the presented approach is that it provides a unified framework
for very general distributed optimization problems.}
Additionally, similar techniques can 
be applied to distributed control problems. We also
presented a simplified algorithm for the design of 
approximating inputs that exploits the problem structure.
{Summarizing, the presented approach provides a systematic 
way to address distributed optimization
problems under mild assumptions on the communication graph
as well as the problem structure.}
{We emphasize that, for practical implementations,
there is still a long way to go. In particular, the highly oscillatory
nature of the approximating inputs as well as time synchronization
will be a major challenge. Apart from that,}
the design of suitable approximating inputs with
improved transient and asymptotic behavior is {complex and} still
an {important} issue to be addressed. {While 
filters can be used as a simple remedy to this problem,}
there are also two other ways we plan to 
approach this problem: (1) altering the choice of
admissible vector fields and (2) modifying the design of
the approximating inputs including an optimal choice
of parameters.

\section*{Acknowledgements}                               
\simonResub{We thank Raik Suttner for his very valuable comments.}

\begingroup
\setlength{\parskip}{0pt}
\setlength{\itemsep}{0pt}
\bibliographystyle{IEEEtran}   
\bibliography{subfiles/bibfile}        
\endgroup

\appendix
\section{Appendix}
\ifthenelse{\boolean{longVersion}}{
\subsection{Proof of \Cref{lemmaSaddlePointDynamics}}\label{proofLemmaSaddlePointDynamics}
{The proof follows a similar argument as the one in~\cite[Theorem 5.1.3]{duerr2015thesis}. First, using~\eqref{eqSPAc}, we have}
	\begin{align}
		\lambda_i(t) = \exp\big( {\int_{0}^{t} \big( a_i x(\tau) - b_i \big) \mathrm{d} \tau } \big) \lambda_i(0),
		\label{eqProofLambdaDyn}
	\end{align}
	{for all $ i = 1,2,\dots, n $;}
	{hence,} $ \lambda_i(0) > 0 $ implies that $ \lambda_i(t) > 0 $,
	for all $ t \geq 0 $, and {consequently,} the set $ \mathcal{R}(\mathcal{M}) $ is positively invariant {w.r.t.~\eqref{eqProofLambdaDyn}}.
	Let $ (x^\star,\dualEq^\star,\dualIneq^\star) $ be an arbitrary point in $ \mathcal{M} $.
	Consider the {candidate} Lyapunov function $ V: \real^n \times \real^{n} \times {\realPos} 
	\to {\realNonNeg} $ defined as
	\begin{align}
		V(x,\dualEq,\dualIneq) 
		&= \tfrac{1}{2} \Vert x-x^\star \Vert^2 + \tfrac{1}{2} \Vert \dualEq - \dualEq^\star \Vert^2 \nonumber \\
		&+ \sum\limits_{i=1}^{n} (\dualIneq_{i} - {\dualIneq_{i}}^\star)
		- \sum\limits_{i: \dualIneq_{i}^\star \neq 0}^{n} \dualIneq_{i}^* \ln( \tfrac{\dualIneq_{i}}{\dualIneq_{i}^\star} ).
	\end{align}
	{We first observe} that $V$ is positive definite with respect to $ (x^\star,\dualEq^\star,\dualIneq^\star) $
	on $ \mathcal{R}(\mathcal{M}) $, {and that all the} level sets are compact. To see this,
	note that according to \cite[p. 207, eq. (1.5)]{bregman1967relaxation}, the
	function $ D: {\realPos^n} \times {\realPos^n} \to \real $
	defined as
	\begin{align}
		D(\dualIneq^\star,\dualIneq)
		&= \sum\limits_{i=1}^{n} \big( \dualIneq_i - \dualIneq_i^\star + \dualIneq_i^\star ( \ln(\dualIneq_i^\star) - \ln(\dualIneq_i) ) \big) \\
		&= \sum\limits_{i=1}^{n} (\dualIneq_{i} - \dualIneq_{i}^\star)
		- \sum\limits_{i: \dualIneq_{i}^\star \neq 0}^{n} \dualIneq^\star_i \ln( \tfrac{\dualIneq_{i}}{\dualIneq_{i}^\star} )
	\end{align}
	is positive for all $ (\dualIneq^\star,\dualIneq) \in {\realPos^n} \times {\realPos^n}  $
	and zero if and only if $ \dualIneq = \dualIneq^* $ \cite[Condition I.]{bregman1967relaxation}
	and its level sets are compact \cite[Condition V.]{bregman1967relaxation}. Thus,
	with $ V(x,\dualEq,\dualIneq) $ additionally being quadratic in $ x $ and $ \dualEq $,
	positive definiteness and compactness of all level sets follows {and hence}
	$V$ is uniformly unbounded on $ \mathcal{R}(\mathcal{M}) $.
	The derivative of $V$ along the trajectories of \eqref{eqSPA} is then given by
	\begin{align}
		&\dot{V}(x,\dualEq,\dualIneq) \\
		\ifthenelse{\boolean{longVersion}}{
		&= -(x-x^\star)^\top \big( \nabla F(x) + \eqConstrL^\top \dualEq + \ineqConstrL^\top \dualIneq \big) \nonumber \\
		&+ (\dualEq-\dualEq^\star)^\top( \eqConstrL x-\eqConstrR)  \nonumber \\
		&- \sum\limits_{ \substack{ i = 1 \\ i \notin \setEq} }^{n} \nu_i^2 + \sum\limits_{i=1}^{n} \dualIneq_{i} ( \ineqConstrLAgent_i x - \ineqConstrR_i ) 
		- \sum\limits_{i: \dualIneq_{i}^\star \neq 0}^{n} 
		\dualIneq_{i}^\star ( \ineqConstrLAgent_i x - \ineqConstrR_i ) \nonumber \\}{}
		&= -(x-x^\star)^\top \nabla F(x) - \dualEq^\top \big( \eqConstrL x- \eqConstrR - (\eqConstrL x^\star-\eqConstrR) \big) \nonumber \\
		&- \dualIneq^\top \big( { \ineqConstrL x - \ineqConstrR - (\ineqConstrL x^\star - \ineqConstrR ) } \big) 
		+ (\dualEq-\dualEq^\star)^\top( \eqConstrL x-\eqConstrR) \nonumber \\
		&- \sum\limits_{ \substack{ i = 1 \\ i \notin \setEq} }^{n} \nu_i^2 
		+ \sum\limits_{i=1}^{n} (\dualIneq_{i}-\dualIneq_{i}^\star) ( \ineqConstrLAgent_i x - \ineqConstrR_i )  + F(x) - F(x), \nonumber
	\end{align}
	Using strict
	convexity of $F$, we now have that $ -(x-x^\star)^\top \nabla F(x) < F(x^\star) - F(x) $,
	for all $ x \neq x^\star $ {and hence} 
	we obtain for all $ x \neq x^\star $
	\ifthenelse{\boolean{longVersion}}{
	\begin{align}
		&\dot{V}(x,\dualEq,\dualIneq) \\
		&< F(x^\star)-F(x) - \dualEq^\top \big( \eqConstrL x-\eqConstrR - (\eqConstrL x^\star-\eqConstrR) \big) \nonumber \\
		&- \dualIneq^\top \big( \ineqConstrL x - \ineqConstrR -( \ineqConstrL x^\star- \ineqConstrR) \big) 
		+(\dualEq-\dualEq^\star)^\top( \eqConstrL x-\eqConstrR) \nonumber \\
		&- \sum\limits_{ \substack{ i = 1 \\ i \notin \setEq} }^{n} \nu_i^2 
		+ \sum\limits_{i=1}^{n} (\dualIneq_{i}-{\dualIneq_{i}}^\star) ( \ineqConstrLAgent_i x - \ineqConstrR_i ) + F(x) - F(x) \nonumber \\
		&= \lagrangian(x^\star,\dualEq,\dualIneq) - \lagrangian(x,\dualEq,\dualIneq) + \lagrangian(x,\dualEq,\dualIneq) - \lagrangian(x,\dualEq^\star,\dualIneq^\star) \nonumber \\
		& - \sum\limits_{ \substack{ i = 1 \\ i \notin \setEq} }^{n} \nu_i^2  \nonumber \\
		&= \lagrangian(x^\star,\dualEq,\dualIneq) - \lagrangian(x,\dualEq^\star,\dualIneq^\star) - \sum\limits_{ \substack{ i = 1 \\ i \notin \setEq} }^{n} \nu_i^2.
	\end{align}
	}
	{
	\begin{flalign}
		\dot{V}(x,\dualEq,\dualIneq) &< \lagrangian(x^\star,\dualEq,\dualIneq) - \lagrangian(x,\dualEq^\star,\dualIneq^\star) - \sum\limits_{ \substack{ i = 1 \\ i \notin \setEq} }^{n} \nu_i^2.\hspace*{-1em} & 
	\end{flalign}
	}
	Due to the saddle point property \eqref{eqSaddlePointProperty} {the derivative of $ V $ along the flow} 
	is {strictly negative}, for all $ (x,\dualEq,\dualIneq) $
	except for $ (x,\dualEq,\dualIneq) \in \mathcal{M} $; thus,
	$ (x^\star,\dualEq^\star,\dualIneq^\star) $ is stable 
	according to \cite[Theorem 2.2.2]{duerr2015thesis}. This 
	procedure can be repeated for any point $ (x^\star,\dualEq^\star,\dualIneq^\star) \in \mathcal{M} $,
	hence $ \mathcal{M} $ is stable. Let $ \lagrangian_{\textup{orig}} $
	denote the Lagrangian associated to the original problem
	\eqref{eqOptimizationProblem} and let $ \mathcal{S}_{\textup{orig}} $
	denote the corresponding set of saddle points. Observe
	that $ \lagrangian(x,\dualEq,\dualIneq) = \lagrangian_{\textup{orig}}(x,\dualEq_{\setEq},\dualIneq_{\setIneq}) - \sum_{ { i=1, i \notin \setIneq } }^{n} \dualIneq_i \ineqConstrR_i $
	{ such that $ \dualIneq_i^\star = 0 $ for all $ i = 1,2,\dots,n $, $ i \notin \setIneq $,
	for any saddle point $ (x^\star,\dualEq^\star,\dualIneq^\star) $
	of $ \lagrangian $, since {$ \ineqConstrR_i > 0 $} for $ i = 1,2,\dots,n $, $ i \notin \setIneq $.}
	Thus, the set of saddle points of $ \lagrangian $ is
	given by 
	\begin{align}
		\mathcal{S} = \lbrace~& (x,\dualEq,\dualIneq) \in \real^n \times \real^n \times {\realNonNeg^n} : \\
		&(x,\dualEq_{\setEq},\dualIneq_{\setIneq}) \in \mathcal{S}_{\textup{orig}}, 
		\dualIneq_i = 0 \text{ for } i \notin \setIneq \rbrace \nonumber
	\end{align}
	and hence,
	$\mathcal{M} = \lbrace (x,\dualEq,\dualIneq) \in \real^n \times \real^n \times \real^n : (x,\dualEq_{\setEq},\dualIneq_{\setIneq}) \in \mathcal{S} \text{ and } \dualEq_i = 0 \text{ for } i \notin \setEq, \dualIneq_i = 0 \text{ for } i \notin \setIneq \rbrace $.
	Since $ \mathcal{S}_{\textup{orig}} $ is compact due to \Cref{assMFCQ},
	the set $ \mathcal{M} $ is compact as well. {The same argument as the one in the proof of~\cite[Theorem 5.1.3]{duerr2015thesis}
	then yields that the} set of saddle points is 
	asymptotically stable with respect to the set of initial conditions
	$ \mathcal{R}(\mathcal{M}) $. 
}{}

\subsection{Proof of \Cref{lemmaConstructionBrackets}}\label{secAppendixProofConstructionBrackets}
We prove the result by induction. For paths of the form
$ p_{i_1 i_2} = \pathLeft v_{i_1} \pathSep v_{i_2} \pathRight $, i.e., $ \length(p_{i_1 i_2}) = 1 $, by 
\eqref{eqDefAdmissibleVectorFields} and \eqref{eqDefRecursionBStart}  
equation \eqref{eqLemmaBracket} follows immediately.
Further we observe that the vector field 
\eqref{eqDefRecursionBStart} is admissible if $ k_1 \in \mathcal{I}(j), k_2 \in \mathcal{I}(i) $
and $ g_{ij} \neq 0 $, which is true since $ p_{ij} $ is a path in $ \mathcal{G} $.
Suppose now that the result holds
for all paths $ p $ with $ \length(p) \leq \bar{\ell} $,
$ \bar{\ell} \geq 2 $. Let $p_{i_{1} i_{\ell}} = \pathLeft v_{i_{1}} \pathSep v_{i_{2}} \pathSep \dots \pathSep v_{i_{k}} \pathRight $
be any path with $ \length( p_{i_1 i_k} ) = \bar{\ell} + 1 $.
Let further $ q_r \in \subpath_{i_{1} \bullet}(p_{i_{1} i_{k}}) $
be a subpath of $p_{i_{1} i_{k}} $ that ends at $ v_r $, $ r = i_2, i_3, \dots, i_{k-1} $.
Then, since $ \length(q_r)\leq \bar{\ell} $, $ \length(q_r^c)\leq\bar{\ell} $,
we have by \eqref{eqDefRecursionBStep} and the induction hypothesis
\begin{align}
    \recBracket{k_1,k_2}\big(p_{i_1 i_k}\big)(\state) 
    &=
    \big[ \recBracket{k_1,s}{}(q_r^c), \recBracket{s,k_2}{}(q_r) \big](\state) \nonumber \\
    &= 
    \ifthenelse{\boolean{longVersion}}{
    \big[ h_{k_1,s}, h_{s,k_2} \big](\state) \nonumber \\
    &=
    \unitVec{k_2} \unitVec{s}^\top \state_{k_1} \unitVec{s} - \unitVec{s} \unitVec{k_1}^\top \state_{s} \unitVec{k_2} \nonumber \\
    &=
    \state_{k_1} \unitVec{k_2},}{
    \big[ h_{k_1,s}, h_{s,k_2} \big](\state) = 
    \state_{k_1} \unitVec{k_2}, 
    }
\end{align}
where $ s \in \mathcal{I}(\tail(q_r)) $ and where we {have} used that $ k_1 \neq k_2 $.
This proves \eqref{eqLemmaBracket}. Further, if $ k_1 \in \mathcal{I}\big( \tail(p_{i_1 i_k}) \big) $, 
i.e., $ k_1 \in \mathcal{I}\big( \tail(q_r^c) \big) $, then,
by the induction hypothesis and with $ s \in \mathcal{I}(\tail(q_r)) = \mathcal{I}(\head(q_r^c))$,
$ \recBracket{k_1,s}{}(q_r^c) $ is a Lie bracket
of admissible vector fields. Similarly, if $ k_2 \in \mathcal{I}\big( \head(p_{i_1 i_k}) \big) $,
by the induction hypothesis and with $ s \in \mathcal{I}(\tail(q_r)) $,
also $ \recBracket{s,k_2}{}(q_r) $ is a Lie bracket of admissible 
vector fields. Thus, $ \recBracket{k_1,k_2}\big(p_{i_1 i_k}\big) $
is a Lie bracket of admissible vector fields as well{,} which
concludes the proof.

\subsection{Proof of \Cref{lemmaChoiceSubpath}}\label{secAppendixProofChoiceSubpath}
{We first observe first} that \eqref{eqRecPHall} is the same as \eqref{eqDefRecursionBStart},
	\eqref{eqDefRecursionBStep} with a special choice of the subpath as well
	as an additional projection with the property $ \resort_{{\mathbb{P}}}\big( B \big)(\state) = B(\state) $
	for all $ \state \in \real^{3n} $. {Hence, it immediately follows that 
	$ \recBracketPHall{k_1,k_2}{}(p_{i_1i_r})(\state) = \recBracket{k_1,k_2}{}(p_{i_1i_r})(\state) $}. In the same manner, {we also have that}
	\ifthenelse{\boolean{longVersion}}{
	\begin{align}
        \degree\big( \recBracketPHall{k_1,k_2}{}(p_{i_1i_r}) \big) = \degree\big( \recBracket{k_1,k_2}{}(p_{i_1i_r}) \big) = \length(p_{i_1i_r}).
	\end{align}
	}{
	{$
		\degree\big( \recBracketPHall{k_1,k_2}{}(p_{i_1i_r}) \big) = \degree\big( \recBracket{k_1,k_2}{}(p_{i_1i_r}) \big) = \length(p_{i_1i_r}).
	$}}
	We show {the second part} by induction. First observe that
	for paths $ p_{i_1 i_r} $ with $ \length(p_{i_1 i_r}) = 1 $ it is 
	clear that $ \recBracketPHall{k_1,k_2}{}(p_{i_1 i_r}) \in \mathbb{P} $
	since $ \recBracket{k_1,k_2}{}(p_{i_1 i_r}) $ is an admissible vector field
	by \Cref{lemmaConstructionBrackets} and all admissible vector fields are in $ \mathbb{P} $.
	{Further, for paths $ p_{i_1 i_r} $ with $ \length(p_{i_1 i_r}) \in \lbrace 2,3,4,6 \rbrace $} 
	it also {follows from the definition of the projection} that $ \recBracketPHall{k_1,k_2}{}(p_{i_1 i_r}) \in \mathbb{P} $. Suppose now that the result
	holds true for all paths $p$ with $ \length(p) = \bar{\ell} $, where $ \bar{\ell} \geq 2 $,
	and consider a path $ p_{i_1 i_r} $ with $ \length(p_{i_1 i_r}) = \bar{\ell} + 1 $.
	Observe that all subbrackets of $ \recBracketPHall{k_1,k_2}{}(p_{i_1 i_r}) $ are in $ \mathbb{P} $
	by the induction hypothesis and hence, by \ref{itemPH3}, \ref{itemPH4a}, \ref{itemPH4b}, we have $ {B} := \recBracketPHall{k_1,k_2}{}(p_{i_1i_r}) \in \mathbb{P} $
	if 
	\begin{align}
		 \degree\big( \lFac( {B} ) \big) &< \degree\big( \rFac( {B} ) \big) \label{eqProofPHallCond1} \\
		 \degree\big( \lFac( \rFac( {B} ) ) \big) &<  \degree\big( \lFac( {B} ) \big) \label{eqProofPHallCond2}; 
	\end{align}
	{we will} show next that these conditions are fulfilled for the above choice
	of subpaths.	
	By \eqref{eqRecPHall} and \eqref{eqChoiceSubpath} we have  
	\begin{align*}
		\degree\big( \rFac( \recBracketPHall{k_1,k_2}{}(p_{i_1i_r}) ) \big) &= \degree\big( \recBracketPHall{s,k_2}{}(q) \big) = \length(q) \\
		\degree\big( \lFac( \recBracketPHall{k_1,k_2}{}(p_{i_1i_r}) ) \big) &= \degree\big( \recBracketPHall{k_1,s}{}(q^c) \big) = \length(p_{i_1 i_r}) - \length(q).
	\end{align*}
	{Since} $ \lfloor \tfrac{a}{b} \rfloor \geq \tfrac{a-b+1}{b} $, for all $ a \in \mathbb{Z}, b \in \nat $,
	we infer
	\ifthenelse{\boolean{longVersion}}{
	\begin{align}
		\length(q) &= \theta(p_{i_1 i_r}) -1 \geq \tfrac{\length(p_{i_1 i_r}) + 1}{2},
	\end{align}}{
	{$ \length(q) = \theta(p_{i_1 i_r}) -1 \geq \tfrac{1}{2}(\length(p_{i_1 i_r}) + 1)$},}
	for $ \length(p_{i_1 i_r}) \geq 5 $, and hence we obtain 
	\begin{align}
		\degree\big( \rFac( \recBracketPHall{k_1,k_2}{}(p_{i_1i_r}) ) \big) - \degree\big( \lFac( \recBracketPHall{k_1,k_2}{}(p_{i_1i_r}) ) \big) \nonumber \\
		\geq \length(p_{i_1 i_r}) + 1 - \length(p_{i_1 i_r}) > 0.
	\end{align}
	Thus, \eqref{eqProofPHallCond1} holds. For \eqref{eqProofPHallCond2}, we
	first note that
	\ifthenelse{\boolean{longVersion}}{
	\begin{flalign}
		\degree\big( \lFac( \rFac( \recBracketPHall{k_1,k_2}{}(p_{i_1i_r}) ) ) \big)
		&= \degree\big( \lFac( \recBracketPHall{s,k_2}{}(q) ) \big) \hspace*{-1em} & 
	\end{flalign}
	}
	{$ \degree\big( \lFac( \rFac( \recBracketPHall{k_1,k_2}{}(p_{i_1i_r}) ) ) \big) = \degree\big( \lFac( \recBracketPHall{s,k_2}{}(q) ) \big) $}
	and, since $ \lFac( \recBracketPHall{s,k_2}{}(q) ) \in \mathbb{P} $ by the
	induction hypothesis, it is 
	$ \degree\big( \lFac( \recBracketPHall{s,k_2}{}(q) ) \big) \leq \degree\big( \rFac( \recBracketPHall{s,k_2}{}(q) ) \big) = \length(q) - \degree\big( \lFac( \recBracketPHall{s,k_2}{}(q) ) \big) $
	according to \ref{itemPH4a}. Hence, we obtain
	\ifthenelse{\boolean{longVersion}}{
	\begin{align}
		\degree\big( \lFac( \rFac( \recBracketPHall{k_1,k_2}{}(p_{i_1i_r}) ) ) \big) \leq \tfrac{\length(q)}{2}.
	\end{align}
	}{$ \degree\big( \lFac( \rFac( \recBracketPHall{k_1,k_2}{}(p_{i_1i_r}) ) ) \big) \leq \tfrac{\length(q)}{2}. $}
	{As a result},~\eqref{eqProofPHallCond2} is fulfilled {when}
	\begin{align}
		\tfrac{\length(q)}{2} \leq \length(p_{i_1 i_r}) - \length(q). \label{eqProofPHallCond2Reformulated}
	\end{align}
	{We now compute}
	\begin{align*}
		\tfrac{3}{2} \length(q) &= \tfrac{3}{2} \lfloor \tfrac{\length(p_{i_1 i_r})}{2} \rfloor +\tfrac{3}{2} \leq \tfrac{3}{4} \length(p_{i_1 i_r}) + \tfrac{3}{2} \leq \length(p_{i_1 i_r}),  
	\end{align*}
	for $ \length(p_{i_1 i_r}) \geq 6 $; for $ \length(p_{i_1 i_r}) = 5 $, 
	we have that $ \tfrac{3}{2} \length(q) = \tfrac{9}{2} <  \length(p_{i_1 i_r}) $, thus
	\eqref{eqProofPHallCond2Reformulated} holds for all considered $ p_{i_1 i_r} $
	which proves that \eqref{eqProofPHallCond2} holds; {this concludes the proof.}

\subsection{Proof of \Cref{lemmaLowDegreeBrackets}}\label{secProofLemmaLowDegreeBrackets}
It is clear that \eqref{eqLemmaLowDegree} holds for $ \length(p_{ij}) = 2 $, since
	$ \recBracketPHall{r+j,i}{}( p_{ij} ) $ is a bracket of degree two, i.e., a bracket
	of the form $ [ \phi_{k_1}, \phi_{k_2} ] $, $ k_1 \neq k_2 $, such that
	\begin{align}
		{E_{\resort_{{\mathbb{P}}}( [ \phi_{k_1}, \phi_{k_2} ] ) }}
		&=
		\begin{cases}
			[ \phi_{k_1}, \phi_{k_2} ] \text{ if } k_1 < k_2 \\
			[ \phi_{k_2}, \phi_{k_1} ] \text{ if } k_2 < k_1.
		\end{cases}  
	\end{align}
	Consider now a path $ p_{i_1 i_4} = \pathLeft v_{i_1} \pathSep v_{i_2} \pathSep v_{i_3} \pathSep v_{i_4} \pathRight $, $ i_1 \neq i_2 \neq i_3 \neq i_4 $,
	i.e., $ \length(p_{i_1 i_4}) = 3 $. Then
	\begin{align}
		&\recBracketPHall{r+i_4,i_1}{}( p_{i_1 i_4} ) \nonumber \\
		&= \resort_{{\mathbb{P}}}\big( \recBracket{r+i_4,i_1}{}( p_{i_1 i_4} ) \big) \nonumber \\
		&= \resort_{{\mathbb{P}}}\big( \big[ h_{r+i_4,r+i_3}, [ h_{r+i_3,r+i_2}, h_{r+i_2,i_1} ] \big] \big) \nonumber \\
		&= - \big[ h_{r+i_4,r+i_3}, [  h_{r+i_2,i_1}, h_{r+i_3,r+i_2} ] \big],
	\end{align}
	where we {have} used the assumption on the ordering of the P.~Hall basis.
	The only equivalent bracket in $ \mathbb{P} $ is then given by
	$ {B} = \big[ h_{r+i_3,r+i_2} , [  h_{r+i_2,i_1}, h_{r+i_4,r+i_3} ] \big] $,
	but {we have that} $ B(\state) \equiv 0 $, since
	\begin{align}
		&[  h_{r+i_2,i_1}, h_{r+i_4,r+i_3} ](\state) \nonumber \\
		&= \unitVec{r+i_3} \unitVec{r+i_4}^\top \unitVec{i_1} \state_{r+i_2} - \unitVec{i_1} \unitVec{r+i_2}^\top \unitVec{r+i_3} \state_{+i_2}
		= 0.
	\end{align}
	Thus, the claim follows.

\subsection{A simplified algorithm for the construction of approximating sequences}\label{secAlgorithmApproxSeq}
\ifthenelse{\boolean{longVersion}}{
\begin{table*}[ht]
	\begin{center}
		\begin{tabular}{ccccccccccccccccc}
			\toprule 
			$ \degree(B) $            & $ 2 $ & $ 3 $ & $ 4 $ & $ 5 $ & $ 6 $ & $ 7 $ & $ 8 $ & $ 9 $ & $ 10 $ & $ 11 $ & $ 12 $ & $13$ & $14$ & $15$ & $16$ & $17$\\
			\midrule
			$\vert {E_{B,\textup{full}}} \vert$        & $1$ & $2$ & $3!$ & $4!$ & $5!$ & $6!$ & $7!$ & $8!$ & $9!$ & $10!$ & $11!$ & $12!$ & $13!$ & $14!$  & $15!$  & $16!$ \\
			$\vert {E_B} \vert $ & $1$ & $1$ & $2$ & $3$  & $5$   & $9$   & $16$   & $28$    & $51$ & $93$  & $170$ & $315$ & $585$ & $1089$ & $2048$ & $3855$\\
			\bottomrule    
		\end{tabular}
	\end{center}
	\caption{A comparison of $ \vert {E_{B,\textup{full}}} \vert $ and $ \vert {E_B} \vert $ for a specific choice
	of the P.~Hall basis that fulfills the assumptions as in \Cref{lemmaLowDegreeBrackets}.
	The numbers were obtained 	by symbolically computing the resulting vector fields using a computer algebra
	system. Interestingly, the sequence of $ \vert {E_B} \vert $ has two
	matching sequences \cite{oeis2017series1} and \cite{oeis2017series2}
	except for the value for $ \degree(B) = 15 $ which should be
	1091 or 1092, thus we conjecture that these
	sequences are a good upper bound for $ \vert {E_B} \vert $.	
	} \label{tableNonZeroBrackets}
\end{table*}}{}
{Our objective in this section is to provide} a modified version of the construction
procedure from~\cite{liu1997approximation} {using the structural
properties of the problem at hand, which leads to considerable simplifications.}
{Given the scopes of this paper and the complicated nature of the subject,} we do not discuss this algorithm in detail; we refer the reader to~\ifthenelse{\boolean{longVersion}}{the supplementary material of the present manuscript}{\cite{mic2017distributedArxiv} and the references therein}, as well as the original work~\cite{liu1997approximation}.
We first {provide a formal definition of the already mentioned equivalence relation}
on the set of Lie brackets:
\begin{definition}[Equivalent brackets]\label{defEquivalenceBrackets}
Let $ \mathcal{PH} = ( \mathbb{P}, \prec ) $ be a P.~Hall
	basis of $ \Phi = \lbrace \phi_1, \dots, \phi_M \rbrace $
	and let $ \degree_k(B) $ denote the degree of the vector field $ \phi_k $ in the
	bracket $ B \in \mathcal{PH} $. We say that two brackets 
	$ B_1, B_2 \in {\mathbb{P}} $ are equivalent, {denoted by} $ B_1 \sim B_2 $,
	if $ \degree_k(B_1) = \degree_k(B_2) $ for all $ k = 1,\dots,M $.
	\oprocend
\end{definition}
{For a given set of brackets $ \mathbb{P} $, we then denote
by $ E_B = \lbrace \tilde{B} \in \mathbb{P} : \tilde{B} \sim B \rbrace $
the equivalence class corresponding to the bracket $ B \in \mathbb{P} $.
Note that, by definition of the equivalence relation, all brackets
contained in an equivalence class $ E = \lbrace B_1, B_2, \dots, B_r \rbrace $, $ r \in \natPos $,
have the same degree and we hence let $ \degree(E) = \degree(B_k) $, $ k \in \lbrace 1,2,\dots,r \rbrace $, denote 
the degree of the equivalence class.}
For the construction of the sets of frequencies, we {also need the} following two definitions:
\begin{definition}[Minimally canceling] \label{defSymmetricallyMinimallyCanceling}
A set $ \Omega = \lbrace \omega_1, \dots, \omega_m \rbrace $
	is called minimally canceling if for each collection of integers {$ \lbrace y_i \rbrace_{i=1}^m $},
	such that $ \sum_{k=1}^m \vert y_k \vert \leq m $ we have 
	$ \sum_{k=1}^{m} y_k \omega_k = 0 $ if and only if all $ y_k $ are equal. 
	\oprocend
\end{definition}
\begin{definition}[Independent collection] \label{defIndependentCollection}
A finite collection of sets  {$ \lbrace \Omega_{\lambda} \rbrace_{\lambda=1}^{N} $}, where
	$ \Omega_{\lambda} = \lbrace \omega_{\lambda,1}, \omega_{\lambda,2} \dots, \omega_{\lambda, M_{\lambda}} \rbrace $, is
	called independent if the followings hold: \\[-2em]
	\begin{enumerate}[leftmargin=*]
	\item the sets $ \Omega_{\lambda} $ are pairwise disjoint, and
	\item for each 
	collection of integers {$ \lbrace y_{i,k} \rbrace_{i=1}^{N} $}, $ k~=~1,\dots,M_i $, such that
	\begin{align*}
		\sum\limits_{i=1}^{N} \sum\limits_{k=1}^{M_i} y_{i,k} \omega_{i,k} = 0 \quad \text{and} \quad \sum\limits_{i=1}^{N} \sum\limits_{k=1}^{M_i} \vert y_{i,k} \vert \leq \sum\limits_{i=1}^{N} M_i
	\end{align*}
	{we have}
	\ifthenelse{\boolean{longVersion}}{
	\begin{align}
		\sum\limits_{k=1}^{M_i} y_{i,k} \omega_{i,k} = 0,
	\end{align}}{
	$\sum_{k=1}^{M_i} y_{i,k} \omega_{i,k} = 0$,}
	for each $ i = 1,\dots,N $. \oprocend
	\end{enumerate}
\end{definition}
Consider now an extended system of the form
\begin{align}
	\dot{\state} = f_0(\state) + \sum\limits_{\substack{ B \in \mathcal{B} \\ \degree(B) \geq 2 }} v_B B(\state),
	\label{eqAppendixAlgoExtendedSystem}
\end{align}
where $ f_0: \real^N \to \real^N $, $ \mathcal{B} \subset \mathbb{P} ${, $ \mathcal{B} $ finite,} for some
P.~Hall basis $ \mathcal{PH}(\Phi) = ( \mathbb{P}, \prec ) $,
$ \Phi = \lbrace \phi_1,\phi_2,\dots,\phi_M \rbrace $, $ \phi_k: \real^N \to \real^N $,
$ f_0, \phi_k $ sufficiently smooth, $ v_B \in \real \setminus \lbrace 0 \rbrace $
{and $ B(z) \not\equiv 0 $ for all $ B \in \mathcal{B} $}.
Suppose that for any $ B \in \mathcal{B} $, {we have that} $ \degree_k(B) \in \lbrace 0,1 \rbrace $,
$ k = 1,2,\dots,M $. Consider the system
\begin{align}
	\dot{X}^{\seqParam} = f_0(X^{\seqParam}) + \sum\limits_{k=1}^{M} \phi_k(X^{\seqParam}) U_k^{\seqParam}(t).
	\label{eqAppendixAlgoOriginalSystem}
\end{align}
The following algorithm allows to compute suitable 
input functions $ U_k^{\seqParam} $ such that the solutions
of \eqref{eqAppendixAlgoOriginalSystem} uniformly
converge to those of \eqref{eqAppendixAlgoExtendedSystem}
with increasing $ \seqParam $. It should as well be mentioned
that we also provide an exemplary implementation of the
algorithm in Matlab {in the supplementary material\ifthenelse{\boolean{longVersion}}{}{~of~\cite{mic2017distributedArxiv}}.}
\subsubsection*{Algorithm} 
\textbf{Step 1 {(Determining the equivalence classes)}:} For all $B \in \mathcal{B} $, determine the associated (reduced)
	equivalence class
	\begin{align*}
		{E_B} &= \lbrace \tilde{B} \in {\mathbb{P}} : \tilde{B} \sim B, \tilde{B}(z) \not\equiv 0 \rbrace \\
		&= \lbrace \tilde{B}_{E,1}, \tilde{B}_{E,2}, \dots, \tilde{B}_{E, \vert E(B) \vert} \rbrace,
	\end{align*}
	and let $ \mathcal{E} = \lbrace {E_B}, B \in \mathcal{B} \rbrace $. 
	For each {$ B \in \mathbb{P} $}, set
	\begin{align*}
		\tilde{v}_{B}
		=
		\begin{cases}
			v_B & \text{if } B \in \mathcal{B} \\
			0   & \text{otherwise.}
		\end{cases}
	\end{align*}
\textbf{Step 2 {(Determining the frequencies)}:} For all {$E\in\mathcal{E}_2 := \lbrace E \in \mathcal{E}: \degree(E) = 2 \rbrace $}, 
	choose $ \vert \mathcal{E}_2 \vert $ distinct
	frequencies $ \omega_E \in \real \setminus \lbrace 0 \rbrace $, 
	and for all $E\in\mathcal{E}$, $ \degree(E)\geq 3$ 
	choose $ {M} \vert E \vert $ sets
	\begin{align*}
	    \Omega_{E,\rho,k}^+ &= 
	    \begin{cases}
	        { \lbrace \omega_{E,\rho,k} \rbrace}  & \text{if } \degree_k(E) = 1 \\
	        \emptyset          & \text{if } \degree_k(E) = 0  
	    \end{cases} \\
	    \Omega_{E,\rho,k}^- &= - \Omega_{{E},\rho,k}^+,
	\end{align*}
	{$\omega_{E,\rho,k} \in \real \setminus \lbrace 0 \rbrace $,} $ k=1,\dots,M $, $ \rho = 1,\dots, {\vert E \vert} $, 
	such that\\[-2em]
	\begin{enumerate}[leftmargin=*]
		\item For each $ E \in \mathcal{E} $, $ \degree(E) \geq 3 $,
			and each $ \rho = 1,\dots,\vert E \vert $, the set
			$
				\Omega^+_{E,\rho} = \bigcup_{k=1}^{M} \Omega_{E,\rho,k}^+
			$
			is \hyperref[defSymmetricallyMinimallyCanceling]{minimally canceling}.
		\item The collection of sets \\[-1.5em]
			\begin{align*}
				\Big\lbrace \big\lbrace \omega_E, -\omega_E \big\rbrace_{E \in \mathcal{E}_2}, \big\lbrace \Omega^+_{E,\rho} \cup \Omega^-_{E,\rho} \big\rbrace_{\substack{E \in \mathcal{E},\degree(E)\geq3,\\ \rho=1,\dots,{\vert E \vert}} } \Big\rbrace
			\end{align*}\\[-1em]
			is \hyperref[defIndependentCollection]{independent}.
	\end{enumerate}
\textbf{Step 3 {(Calculating the auxiliary matrix $ \Xi_E $)}:} For all $ E \in \mathcal{E} $ with $ \degree(E) \geq 3 $, compute
	\begin{align*}
		\Xi_{E} =
		\begin{bmatrix}
			\xi_{\tilde{B}_{E,1},1}^+               & \xi_{\tilde{B}_{E,1},2}^+               & \dots  & \xi_{\tilde{B}_{E,1},\vert E \vert}^+ \\
			\xi_{\tilde{B}_{E,2},1}^+               & \xi_{\tilde{B}_{E,2},2}^+               & \dots  & \xi_{\tilde{B}_{E,2},\vert E \vert}^+ \\
			\vdots                              & \vdots                              & \ddots & \vdots \\
			\xi_{\tilde{B}_{E,\vert E \vert},1}^+ & \xi_{\tilde{B}_{E,\vert E \vert},2}^+ & \dots  & \xi_{\tilde{B}_{E,\vert E \vert},\vert E \vert}^+
		\end{bmatrix},
	\end{align*}
	where{, for any $ B \in E $, we let}
	\begin{align*}
		\xi_{B,\rho}^{+} &= 
		\hat{g}_B( \omega_{E,\rho,\theta_B(1)}, \omega_{E,\rho,\theta_B(2)}, \dots, \omega_{E,\rho,\theta_B(\degree(B))} ), 
	\end{align*}
	with $ \theta_B(i) = k $ if the $i$th vector field in $B$ is $ \phi_k $ and where
	 $ \hat{g}_B: \real^{\degree(B)} \to \real $ is defined as follows: \\[-2em]
	\begin{itemize}[leftmargin=*]
		\item If $ \degree(B) = 1 $, then
		$
			\hat{g}_B(\tilde{\omega}_1) = 1.
		$
		\item If $ B = [B_1,B_2] $, then
		\begin{align*}
			&\hat{g}_B(\tilde{\omega}_1,\tilde{\omega}_2,\dots,\tilde{\omega}_{\degree(B)}) = \frac{ \hat{g}_{B_1}( \tilde\omega_{1},\tilde\omega_{2}, \dots, \tilde\omega_{\degree(B_1)} ) }{ \sum_{i=1}^{\degree(B_1)} \tilde\omega_{i}  }  \\
			&\times\;
			\hat{g}_{B_2}( \tilde\omega_{ \degree(B_1) +1}, \tilde\omega_{ \degree(B_1) +2}, \dots, \tilde\omega_{ \degree(B_1) + \degree(B_2)} ).  
		\end{align*}
	\end{itemize}
\textbf{Step 4 {(Calculating the input coefficients)}:} For all $ E \in \mathcal{E} $ with $ \degree(E) = 2 $, i.e., $ E(B) = \lbrace B \rbrace = \big[ \phi_{k_1}, \phi_{k_2} \big] $, set
	\begin{align*}
	\eta_{E,k_1}(\omega_E)
	&= \img \tfrac{1}{\beta_E} \text{sign}(\tilde{v}_B \omega_{E}) \sqrt{ \tfrac{1}{2} \vert \tilde{v}_B \omega_{E} \vert } \\
	\eta_{E,k_2}(\omega_E)
	&= \beta_E \sqrt{ \tfrac{1}{2} \vert \tilde{v}_B \omega_{E} \vert },
	\end{align*}
	where $ \beta_E \neq 0 $. 
	For all $ E \in \mathcal{E} $ with $ \degree(E) \geq 3 $ let\footnote{{We tacitly assume here that $\Xi_{E}$ is invertible. It has been shown in~\cite{liu1997approximation} that there always exists a choice of frequencies such that the corresponding matrix obtained when using ``full'' equivalence classes is invertible; however, it is not clear whether this also holds in the case of reduced equivalence classes where $ \Xi_E $ is a submatrix obtained from the general one by removing several rows and columns. }}
	\begin{align*}
		\begin{bmatrix} \gamma_{E,1} \\ \gamma_{E,2} \\ \smash[t]{\vdots} \\[-0.5em] \gamma_{E,\vert E \vert} \end{bmatrix}
		=
		\Xi_{E}^{-1} 
		\begin{bmatrix} \tilde{v}_{\tilde{B}_{E,1}} \\ \tilde{v}_{\tilde{B}_{E,2}}  \\ \smash[t]{\vdots} \\[-0.5em] \tilde{v}_{\tilde{B}_{E,\vert E \vert}}  \end{bmatrix}
	\end{align*}
	and compute $ \eta_{E}(\omega) $ as follows: \\[-2em]
	\begin{itemize}[leftmargin=*]
		\item If $ \degree(E) $ is odd, {for each $\rho = 1,\dots,\vert E \vert$,} take\\[-1.5em]
		\begin{align*}
		\eta_E(\omega) = \beta_{E,\omega} \big( \tfrac{1}{2} \gamma_{E,\rho} \img^{\degree(E)-1} \big)^{\tfrac{1}{\degree(E)}}
		\end{align*}
		for all $ \omega \in \Omega_{E,\rho}^+ $, and
		\item if $ \degree(E) $ is even, {for each $\rho = 1,\dots,\vert E \vert$,} take\\[-1.5em]
		\begin{align*}
		\eta_{E}(\tilde{\omega}) &= \img \beta_{E,{\tilde{\omega}}} \text{sign}(\gamma_{E,\rho}(t) \img^{\degree(E)-2}) \left\vert \tfrac{1}{2} \gamma_{E,\rho}(t) \img^{\degree(E)-2} \right\vert^{\tfrac{1}{\degree(E)}} 
		\end{align*}
		for some $ \tilde{\omega} \in \Omega_{E,\rho}^+ $ and\\[-1.5em]
		\begin{align*}
		\eta_{E}(\omega) &= \beta_{E,\omega} \left\vert \tfrac{1}{2} \gamma_{E,\rho}(t) \img^{\degree(E)-2} \right\vert^{\tfrac{1}{\degree(E)}} 
		\end{align*}
		for all $ \omega \in \Omega_{E,\rho}^+ \setminus \lbrace \tilde{\omega} \rbrace $.
	\end{itemize}
	In both cases $ \beta_{E,\omega} \in \real $ can be chosen freely such that it
	fulfills $ \prod_{\omega \in \Omega_{E,\rho}^+} \beta_{E,\omega} = 1 $. \\[0.5em]
\textbf{Step 5 {(Calculating the approximating inputs)}:} Compute the input according to $ U_k^{\sigma}(t) = \sum_{E \in \mathcal{E}} U_{k,E}^{\sigma}(t) $
	with $ U_{k,E}^{\sigma} : \real \to \real $ being defined as follows: \\[-2em]
	\begin{itemize}[leftmargin=*]
		\item If $ \degree_k(E) = 0 $: $ U_{k,E}^{\sigma}(t) = 0 $.
		\item If $ \degree(E) = 2, \degree_k(E) = 1 $:
		\begin{align*}
			U_{k,E}^{\sigma}(t) = 2 \sqrt{\sigma} \Real \big( \eta_{E,k}(\omega_E) e^{\img \sigma \omega_{E} t} \big) .
		\end{align*}
		\item If $ \degree(E) = N $, $ \degree_k(E) = 1 $:
		\begin{align*}
			U_{k,E}^{\sigma}(t) = 2 \sigma^{\tfrac{N-1}{N}} \sum\limits_{\rho=1}^{\vert E \vert} \Real \big( \eta_{E}(\omega_{E,\rho,k}) e^{\img \sigma \omega t} \big).
		\end{align*}
	\end{itemize}

Note that this algorithm is a reformulation of the one presented
in \cite{liu1997approximation} \ifthenelse{\boolean{longVersion}}{(see the supplements of this manuscript for a derivation)}{(see the supplementary material of~\cite{mic2017distributedArxiv} for a derivation)} exploiting two structural
properties of the problem at hand: (1) each $ B \in \mathcal{B} $ 
fulfills $ \degree_k(B) \in \lbrace 0,1 \rbrace $ for all $ k =1,2,\dots,M $
and (2) {a large number of} the equivalent brackets evaluate to zero \ifthenelse{\boolean{longVersion}}{{(see~\Cref{tableNonZeroBrackets})}}{{see~\cite[Table~2 and Lemma~4]{mic2017distributedArxiv}}}. 
{Note that}~(1) simplifies the calculation of $ \xi^+_{B,\rho} $ in 
step~3 and~(2) reduces the cardinality of each $ {E_B} $ in step~1, where
usually the full equivalence class $ {E_{B,\textup{full}}} = \lbrace \tilde{B} \in \mathcal{B}: \tilde{B} \sim B \rbrace $
is used, thus leading to a reduction of the dimension of $ \Xi_{E} $ in step~3
and hence also simplifying step~4.
\ifthenelse{\boolean{longVersion}}{
In fact, we can derive the
following result on the equivalent brackets:
\begin{lemma}\label{lemmaEquivalentBrackets}
Consider a graph $ \mathcal{G} = ( \mathcal{V}, \mathcal{E} ) $ of $n$ nodes.  
	Let $p_{i_1 i_{r}} = ( v_{i_1}, v_{i_2}, \dots, v_{i_{r}} ) $
	be the shortest path between $v_{i_1}$ and {$v_{i_{r}}$},
	$ v_{i_k} \in \mathcal{V} $ for $ k = 1,2,\dots,r $, $r\geq 3$.
	Let $ \Phi = \lbrace \phi_{a_1}, \phi_{a_2}, \dots, \phi_{a_{r-1}} \rbrace $ 
	be a set of vector fields with
	\begin{align}
		\phi_{a_j} \in \lbrace h_{k_1,k_2}: k_1 \in \mathcal{I}(i_{j+1}), k_2 \in \mathcal{I}(i_{j}) \rbrace,
	\end{align}
	for $ j = 1,2,\dots, r - 1 $.
	Denote some given P.~Hall basis of $ \Phi $ by $ \mathcal{PH}(\Phi) = ( \mathbb{P}, \prec ) $.
	Let $ B \in \mathbb{P} $ and suppose that $ \degree_{a_j}(B) \in \lbrace 0, 1 \rbrace $
	for $ j = 1,2,\dots,r-1 $. Define $ \mathcal{J}(B) = \lbrace j=1,2,\dots,r-1 : \degree_{a_j}(B) = 1 \rbrace $
	and further denote
	\begin{align}
		j_{\textup{min}}(B) = \min\limits_{j \in \mathcal{J}(B)} \lbrace j \rbrace,  \quad j_{\textup{max}}(B) = {\max\limits_{j \in \mathcal{J}(B)} \lbrace j \rbrace} .
	\end{align}
	{
	Then, if $ \mathcal{J}(B) $ is a connected set, i.e., $ \mathcal{J}(B) = \lbrace j_{\textup{min}}(B), j_{\textup{min}}(B)+1, \dots, j_{\textup{min}}(B) + \degree(B) - 1 \rbrace $
	and $ j_{\textup{max}}(B) = j_{\textup{min}}(B) + \degree(B) - 1 $, {for any $ k_1 \in \mathcal{I}(i_{j_{\textup{max}}(B)+1})$, $k_2 \in \mathcal{I}(i_{j_{\textup{min}}(B)}) $ and for all $ z \in \real^{3n} $, we have that} 
	\begin{align}
		B(z) = \pm h_{k_1,k_2}(z) \quad \text{or} \quad B(z) = 0 \label{eqLemmaNonZeroEquivalentBrackets}.
	\end{align}
	If $ \mathcal{J}(B) $ is not a {connected set}, we have $ B(\state) = 0 $
	for all $ z \in \real^{3n} $. 
	}
	\oprocend
\end{lemma}
\begin{proof}
We prove this result by induction. Suppose first that
	$ \degree(B) = 1 $. {Then $ \mathcal{J}(B) = \lbrace j_{\textup{min}} \rbrace = \lbrace j_{\textup{max}} \rbrace $, which means it has only one element.}
	Hence, the claim is obviously true. Since the case of $ \mathcal{J}(B) $
	not being a connected set does not appear for $ \degree(B) = 1 $, we also
	look at $ \degree(B) = 2 $. Let $ \mathcal{J}(B) = \lbrace j_1,j_2 \rbrace $,
	$ j_1 \neq j_2 $. Observe that, for all $ j_1,j_2 = 1,2,\dots,r-1 $, $ j_1 \neq j_2 $, and
	$ j_1 \leq r - 2 $ (or $ j_2 \leq r -1 $), we have 
	\begin{align}
		B(z) 
		&= [\phi_{a_{j_1}}, \phi_{a_{j_2}}](\state) \nonumber \\
		&= [ h_{k_1,k_2}, h_{k_3,k_4} ](\state) \nonumber \\
		&= [ \state_{k_1} \unitVec{k_2}, \state_{k_3} \unitVec{k_4} ] \nonumber \\
		&= \unitVec{k_4} \unitVec{k_3}^\top \unitVec{k_2} \state_{k_1} - \unitVec{k_2} \unitVec{k_1}^\top \unitVec{k_4} \state_{k_3},
	\end{align}
	where $ k_1 \in \mathcal{I}(i_{j_1+1}), k_2 \in \mathcal{I}(i_{j_1}), k_3 \in \mathcal{I}(i_{j_2+1}), $ {and} $ k_4 \in \mathcal{I}(i_{j_2}) $.
	We then compute
	\begin{align}
		[\phi_{a_{j_1}}, \phi_{a_{j_2}}](\state)
		&=
		\begin{cases}
			\state_{k_1} \unitVec{k_4}  & \text{if } k_2 = k_3 \\
			-\state_{k_3} \unitVec{k_2} & \text{if } k_1 = k_4 \\
			0                           & \text{otherwise}.
		\end{cases}
	\end{align}
	{Note that} $ k_2 = k_3 $ only if $ i_{j_1} = i_{j_2 + 1 } $,
	i.e., $ j_1 = j_2 +1 =j_{\textup{max}} $, $ j_{\textup{min}} = j_2 $, and $ k_1 = k_4 $ only if $ i_{j_2} = i_{j_1 + 1} $,
	i.e., $ j_2 = j_1 + 1 = j_{\textup{max}}$, $ j_1 = j_{\textup{min}} $; hence $ B(z) $ is non-zero only if 
	$ \mathcal{J}(B) = \lbrace j_1,j_2 \rbrace $ is connected, 	
	which proves that the claim is true for $ \degree(B) = 2 $.
	Note also that the case $ k_1 = k_4, k_2 = k_3 $ cannot occur since $ j_1 \neq j_2 $. 
	The second claim \eqref{eqLemmaNonZeroEquivalentBrackets} follows immediately
	from these considerations. 
	{To proceed with our induction argument,} suppose now that the claim is true for all $B \in \mathbb{P} $ that
	fulfill the assumptions with $ \degree(B) \leq {\degree}^* $, $ \degree^* \leq r-1 $.
	Consider now some $B \in \mathbb{P} $ with $ \degree(B) = \degree^* + 1 > 2 $.
	Every $B$ can be written as $ B = [B_1,B_2] $,
	where $ \degree(B_1),\degree(B_2) \leq \degree^* $.
	Let $ \mathcal{J}(B) = \lbrace j_1, j_2, \dots, j_{\degree(B)} \rbrace $
	and assume, without loss of generality, that $ j_{k} < j_{k+1} $, for all $ k = 1,\dots,\degree(B)-1 $.
	By the induction hypothesis, $ B_1(\state) $ 	and $ B_2(\state) $ are non-zero only
	if $ \mathcal{J}(B_1) $ and $ \mathcal{J}(B_2) $ are
	both connected sets. Since $ \mathcal{J}(B_2) = \mathcal{J}(B) \setminus \mathcal{J}(B_1) $
	this is the case if and only if 
	{
	\begin{align*}
	\mathcal{J}(B_1)
	& = 
	\begin{cases}
	\lbrace j_1,j_2,\dots,j_{\degree(B_1)} \rbrace  \quad \mathrm{or} \\
	 \lbrace j_{\degree(B)-\degree(B_1)+1},j_{\degree(B)-\degree(B_1)+2},\dots,j_{\degree(B)} \rbrace  
	\end{cases}\\
	& = 
	\begin{cases}
	\lbrace j_1, j_1 + 1, \dots, j_1 + \degree(B_1) - 1 \rbrace, \quad \mathrm{or} \\
	\lbrace j_{\degree(B)-\degree(B_1)+1}, j_{\degree(B)-\degree(B_1)+1} + 1, \dots, \\
	\qquad \qquad \quad j_{\degree(B)-\degree(B_1)+1} + \degree(B_1) - 1 \rbrace. 
	\end{cases}	
	\end{align*}
	}
	We only consider the first case here, since the second case
	can be treated analogously. {Using the first equality above, for $ k_1 \in \mathcal{I}(i_{j_1+\degree(B_1)}), k_2 \in \mathcal{I}(i_{j_1}) $, and 
	$k_3 \in \mathcal{I}(i_{j_{\degree(B)}+1}), k_4 \in \mathcal{I}(i_{j_{\degree(B_1)+1}}) $, we have {by the induction hypothesis} that}
	\begin{align}
		B_1(z) = \pm h_{k_1,k_2}(z) \quad & \text{or} \quad B_1(z) = 0 \\
		B_2(z) = \pm h_{k_3,k_4}(z) \quad & \text{or} \quad B_2(z) = 0.
	\end{align}
	Obviously, following our previous calculations,
	$ [B_1,B_2] $ is non-zero only if $ k_2 = k_3 $, meaning that 
	$ j_1 = j_{\degree(B)}+1$,
	or if $ k_1 = k_4 $, meaning that $ j_1+\degree(B_1) = j_{\degree(B_1)+1} $.
	The first case cannot occur, since $ \degree(B) > 2 $ and $ j_{k+1} > j_k $;
	the second case holds true if and only if $ \mathcal{J}(B) $ is connected,
	thus showing that $ B(z) $ is non-zero only if $ \mathcal{J}(B) $ is connected.
	To show that \eqref{eqLemmaNonZeroEquivalentBrackets} holds,
	consider the case that $ \mathcal{J}(B) $ is connected, i.e.,
	$ \mathcal{J}(B) = \lbrace j_1,j_1 +1, \dots, j_1 + \degree(B) \rbrace $, $ j_{\textup{min}}(B) = j_1 $, $ j_{\textup{max}}(B) = j_1 + \degree(B) $,
	and $ k_1 = k_4 $.
	Then, following the same arguments as before, we have that $ B(z) = \pm h_{k_3,k_2}(z) $
	for $ k_3 \in \mathcal{I}(i_{j_{\degree(B)+1}}) = \mathcal{I}(i_{j_{\textup{max}}(B)+1}) $,
	$ k_2 \in \mathcal{I}(i_{j_1}) = \mathcal{I}(i_{j_{\textup{min}}(B)}) $, which concludes the proof.
\end{proof}
\begin{remark}
The condition that $ \mathcal{J}(B) $ must be a connected set can be interpreted as follows:
	Each admissible vector field $ \phi_{a_j} $ can be associated to an edge in the communication
	graph $ \mathcal{G} $. The condition then means that the vector fields in the bracket
	must be ordered along a path.  
	\oprocend
\end{remark}
}{}

The algorithm presented beforehand still includes several degrees of freedom,
namely the specific choice of frequencies in step 2 as well as the scalings
$ \beta_E, \beta_{E,\omega} $ in step 4. While the conditions on the frequencies
are not hard to satisfy and {in fact, are not restrictive,} it turns out 
that their choice is crucial in practical implementations. {There
is still no constructive way of choosing ``good'' frequencies that we are aware of in the literature}. 
The situation is similar as it comes to the choice of scalings,
but here a heuristic way of how to choose them is to distribute
the energy of the approximating inputs among different admissible 
input vector fields $ \phi_k$. In this spirit, we suggest decreasing the
amplitudes of the approximating inputs entering in the primal variables, 
which will lead to an increase of the amplitudes of the inputs entering
in the dual variables. {Our simulations results indicate that this procedure} 
usually leads to a better transient and asymptotic behavior of the primal variables, which we are typically most
interested in.

\ifthenelse{\boolean{longVersion}}{
{
\subsection{Formal brackets}\label{secAppendixFormalBrackets}
As indicated beforehand, objects such as the degree, the left factor,
the right factor, or a P.~Hall basis are not well-defined for Lie brackets
but need to be defined for formal brackets. We very briefly discuss this
in the following; for a more detailed treatment we refer the reader to 
standard textbooks on the subject, e.g.,~\cite{bourbaki1998lie}, .
Let $ \mathbf{X} = \lbrace X_1, X_2, \dots, X_M \rbrace $ be a finite set of $ M $
non-commuting objects, the so-called indeterminates. We denote by $ \fBr(\mathbf{X}) $
the set of formal brackets constructed from $ \mathbf{X} $, where
a formal bracket is a word fulfilling certain requirements which is constructed from the alphabet 
consisting of the symbols $ X_k $ in $ \mathbf{X} $ as well as the brackets
$ [ $ and $ ] $ and the comma $ , \,$. 
The set of formal brackets $ \fBr(\mathbf{X}) $ is 
then defined as the smallest set of words built from that alphabet 
which contains all elements of $ \mathbf{X} $ and has the property
that, for all $ B_1, B_2 \in \fBr(\mathbf{X}) $, the word $ [ B_1, B_2 ] $
is an element of $ \fBr(\mathbf{X}) $. In this sense, a formal
bracket can be seen as a string representation of a Lie bracket.
However, this string representation is in general not unique.
As an example, we distinguish between the two formal brackets 
$ \big[ \phi_1, [\phi_1,\phi_2] \big] $ and $ \big[ [\phi_2,\phi_1], \phi_1 \big] $,
but these brackets are equivalent as Lie brackets. This is 
the reason why left and right factors as well as the degree
is not well-defined for Lie brackets. For formal brackets
$ B = [ B_1, B_2] \in \fBr( \mathbf{X} ) $, $ B_1, B_2 \in \fBr( \mathbf{X} ) $,
we can uniquely define $ \lFac(B) = B_1 $, $ \rFac(B) = B_2 $ as the left and
right factor of $ B $, respectively. We can further define the degree of 
a formal bracket $ B \in \fBr( \mathbf{X} ) $ in the same way as in~\Cref{secNotation}.

Now, $ \fBr(\mathbf{X}) $ and $ \lieBr(\mathbf{X}) $ are related by 
a mapping $ \mu: \fBr(\mathbf{X}) \to \lieBr(\mathbf{X}) $,
which, in rough words, replaces formal brackets by Lie brackets.
In general, this mapping is not bijective; however, it is if we restrict
the domain of $ \mu $ to a P.~Hall basis of $ \mathbf{X} $~(\cite{bourbaki1998lie}),
which is basically defined in the same way as in~\Cref{defPHallBasis}
but with the set of indeterminates $ \mathbf{X} $ instead of the
the set of vector fields $ \Phi $ and formal brackets instead
of Lie brackets.
Thus, in all of~\Cref{secInputConstruction}, formally we would need to explicitly use $ \mu $ to map from
the formal brackets to Lie brackets as well as an evaluation map
$ \Ev : \lieBr(\mathbf{X}) \to \lieBr(\Phi) $,
which basically simply replaces the indeterminate $ X_i \in \mathbf{X} $
by the vector field $ \phi_i \in \Phi $.
}
}{}

\subsection{Proof of \Cref{theoremMainResult}}\label{secAppendixProofMainResult}
{The proof of \Cref{theoremMainResult} relies on the next general stability result. 
The proof follows the same lines as the proof of~\cite[Theorem~2]{duerr2013Lie}, and is omitted here.}
\begin{lemma}\label{lemmaAppendixPracticalStability}
	Consider {the two dynamics
	\begin{align}
		\dot{\state}(t) &= f\big(t,\state(t)\big), \qquad & \state(t_0) = \state_0,  \label{eqTrajApproxSys1} \\
		\dot{\state}^{\seqParam}(t) &= f^{\seqParam}\big(t,\state^{\seqParam}(t)\big), \qquad & \state^{\seqParam}(t_0) = \state_0, \label{eqTrajApproxSys2}
	\end{align}}
	where $ f, f^{\seqParam}: \real \times \real^n \to \real^n $, 
	$ f,f^{\seqParam} \in \mathcal{C}^1 ${, $ t_0 \in \real $} and $ \seqParam \in {\realPos} $
	is a parameter. Suppose that\\[-2em]
	\begin{enumerate}[leftmargin=*]
		\item a compact set $ \mathcal{S} $ is locally uniformly asymptotically
		stable for \eqref{eqTrajApproxSys1};
		\item the region of attraction $ \mathcal{R}(\mathcal{S}) {~\subseteq \real^n }$ {of $ \mathcal{S} $} is positively
		invariant for \eqref{eqTrajApproxSys2};
		\item for every $ \varepsilon > 0 $, for every $ T > 0 $ and
		{for every $ \mathcal{K}~\subseteq~\mathcal{R}(\mathcal{S}) $}
		there exists $ \seqParam^*>0 $
		such that, for all $ \seqParam > \seqParam^* $, for all $ t_0 \in \real $ {and for all $ z_0 \in \mathcal{K}, $}
		there exist unique solutions $ \state, \state^{\seqParam} $ of
		\eqref{eqTrajApproxSys1} and \eqref{eqTrajApproxSys2} that 
		fulfill 
		\ifthenelse{\boolean{longVersion}}{for all $ t \in  [t_0, t_0 + T] $
		\begin{align}
			\Vert \state(t) - \state^{\seqParam}(t) \Vert \leq \varepsilon.
		\end{align}
		}{$\Vert \state(t) - \state^{\seqParam}(t) \Vert \leq \varepsilon$ for all $ t \in  [t_0, t_0 + T] $.}
	\end{enumerate}
	Then the set $ \mathcal{S} $ is locally practically uniformly
	asymptotically stable for \eqref{eqTrajApproxSys2} and
	$ \state^{\seqParam}(t) $ uniformly converges to $ \state(t) $ on $ [t_0,\infty) $
	for increasing $ \seqParam $.
	\oprocend
\end{lemma}
We are now ready to prove \Cref{theoremMainResult} making use of
\Cref{lemmaAppendixPracticalStability}. 
Since the control law~\eqref{eqDistributedControlInput} is obtained
from the construction procedure presented in~\cite{liu1997approximation},
it follows directly from~\cite[Theorem 8.1]{liu1997approximation}
that for each $ \varepsilon > 0 $, for each $T>0$ and for each initial condition
$ \state^{\sigma}(0) = \state_0 \in \mathcal{R}(\mathcal{M}) $, there exists $ \seqParam^* > 0 $
such that for all {$ \seqParam > \seqParam^* $} and for all $ t \in [0,T] $ 
the inequality \eqref{eqLemmaStabResultTrajApprox} holds{, which shows
convergence on finite time intervals. For the extension to infinite time
intervals we make use of~\Cref{lemmaAppendixPracticalStability}.}
First, note that the set $ \mathcal{M} $ defined by~\eqref{eqDefSetOfSaddlePoints}
is compact by ~\Cref{assMFCQ} (see also the proof of \Cref{lemmaSaddlePointDynamics})
and asymptotically stable for~\eqref{eqSPA} with region of attraction
$ \mathcal{R}(\mathcal{M}) = \lbrace (x,\dualEq,\dualIneq) \in \real^n \times \real^n \times \real^n : \dualIneq \in {\realPos^n} \rbrace $,
according to~\Cref{lemmaSaddlePointDynamics}. 
Also, by the same argumentation as the one in the proof of~\Cref{lemmaSaddlePointDynamics}, the set $ \mathcal{R}(\mathcal{M}) $
is positively invariant for \eqref{eqAllAgentsOriginalSystem}
together with the control law \eqref{eqDistributedControlInput} - \eqref{eqInputBracketHigher}.
For the last assumption in~\Cref{lemmaAppendixPracticalStability},
we first note that we cannot use~\cite[Theorem 8.1]{liu1997approximation},
since, according to assumption 3 in~\Cref{lemmaAppendixPracticalStability}, we are
required to find one $ \seqParam^* $ that works for all 
$ t_0 \in \real $ and for all $ z_0 \in \mathcal{K} $,
but the latter reference only provides uniform convergence in $t$.
However, {for brackets of degree two, by~\cite{duerr2013Lie},
we conclude that such a $ \seqParam^* $ exists;
hence, all assumptions from \Cref{lemmaAppendixPracticalStability}
are fulfilled and the result follows.	
For higher order brackets, the existence of such $\seqParam^*$ has not been shown
explicitly, which is why practical uniform asymptotic stability 
can only be guaranteed if all brackets are of degree two;
still, it is expected that this also holds for the general case.}

\end{document}